\documentclass[10pt,a4paper]{amsart}

\usepackage{amsmath}
\usepackage{amssymb}
\usepackage{amsfonts}
\usepackage{mathrsfs}

\numberwithin{equation}{section}
\newtheorem{theorem}{Theorem}[section]
\newtheorem{corollary}[theorem]{Corollary}
\newtheorem{proposition}[theorem]{Proposition}
\newtheorem{lemma}[theorem]{Lemma}

\theoremstyle{remark}

\theoremstyle{definition}
\newtheorem{definition}{Definition}[section]

\theoremstyle{remark}

\theoremstyle{remark}
\newtheorem{remark}{Remark}[section]

\begin{document}


\title[Rough potential]%
{On the wave equation with a large rough potential}

\author{Piero D'Ancona}
\address{Piero D'Ancona: Unversit\`a di Roma ``La Sapienza'',
Dipartimento di Matematica, Piazzale A.~Moro 2, I-00185 Roma,
Italy}
\email{dancona@mat.uniroma1.it}

\author{Vittoria Pierfelice}
\address{Vittoria Pierfelice: Universit\`a di Pisa, Dipartimento di Matematica,
Via Buonarroti 2, I-56127 Pisa, Italy}
\email{pierfelice@dm.unipi.it}

\begin{abstract}
    We prove an optimal dispersive $L^{\infty}$ decay estimate for a three dimensional wave equation perturbed with a large non smooth potential belonging to a particular Kato class. The proof is based on a spectral representation of the solution and suitable resolvent estimates for the perturbed operator.
\end{abstract}

\thanks{The authors were partially supported
by GNAMPA - Gruppo  "Problemi iperbolici in geometria e fisica"}

\maketitle

\newcommand{\Lui}{\mathscr{L}(L^{1};L^{\infty})}
\newcommand{\Luu}{\mathscr{L}(L^{1};L^{1})}
\newcommand{\Lii}{\mathscr{L}(L^{\infty};L^{\infty})}
\newcommand{\Lpq}{\mathscr{L}(L^{p};L^{q})}
\newcommand{\Lpp}{\mathscr{L}(L^{p};L^{p})}
\newcommand{\Ldd}{\mathscr{L}(L^{2};L^{2})}
\newcommand{\Lz}{L^{\infty}_{0}}
\newcommand{\lapm}{\lambda\pm i\varepsilon}
\newcommand{\lap}{\lambda+i\varepsilon}
\newcommand{\lam}{\lambda-i\varepsilon}
\newcommand{\lapmz}{\lambda\pm i 0}
\newcommand{\lapz}{\lambda+i 0}
\newcommand{\lamz}{\lambda-i 0}
\newcommand{\laep}{\lambda_\varepsilon}
\newcommand{\epsula}{\varepsilon/{2\sqrt{\laep}}}
\newcommand{\wx}{\langle x\rangle}
\newcommand{\wy}{\langle y\rangle}
\newcommand{\supp}{\mathop{\mathrm{supp}}}
\newcommand{\rep}{\rho_{\varepsilon}}
\newcommand{\Vep}{{V_{\varepsilon}}}
\newcommand{\Vth}{{V_{\theta}}}
\newcommand{\Dth}{{D_{V_{\theta}}}}
\newcommand{\Hth}{H_{\theta}}
\newcommand{\uep}{u_{\varepsilon}}
\newcommand{\dist}{\mathop{\mathrm{dist}}}


\section{Introduction}

A basic property of the $n$-dimensional wave equation, $n\geq2$,
\begin{equation}\label{eq.wave}
   \square_{1+n} u=0,\qquad u(0,x)=0,\qquad u_t(0,x)=f(x)
\end{equation}
is expressed by the \emph{dispersive estimate}
\begin{equation}\label{eq.dispfree}
   \|u(t,\cdot)\|_{L^{\infty}}\leq C\;t^{-\frac{n-1}2}
          \|f\|_{\dot B^{\frac{n-1}2}_{1,1}(\mathbb{R}^{n})}.
\end{equation}
Here $\dot B^s_{p,q}(\mathbb{R}^{n})$ is the \emph{homogeneous Besov space}
defined by
\begin{equation}\label{eq.besovfree}
   \|f\|_{\dot B^s_{p,q}(\mathbb{R}^{n})}^{q}
   =\sum_{j\in\mathbb{Z}}2^{jsq}\|\phi_{j}(\sqrt{-\Delta})f\|_{L^{p}}^{q}
\end{equation}
where $\phi_{j}(r)=\phi_{j}(|x|)$ is a Paley-Littlewood partition of unity, i.e., 
$\phi_{j}(r)=\phi_{0}(2^{-j}r)$, $\phi_{0}(r)=\psi(r)-\psi(r/2)$, with $\psi(r)$ 
being a nonnegative function in $C^{\infty}_{0}$ such that $\psi(r)=1$ for $r<1$ 
and $\psi(r)=0$ for $r>2$. 

From \eqref{eq.dispfree} and the energy identity, via interpolation and the 
$T^{*}T$ method, the full set of decay estimates for the wave equation can be 
obtained (see \cite{GV} and 
\cite{KT}). The importance of decay estimates for the applications 
to nonlinear problems is well known.

The possible extension of \eqref{eq.dispfree} to wave equations perturbed with a potential 
\begin{equation}\label{eq.wavepot}
   \square_{1+n} u+V(x)u=0,\qquad u(0,x)=0,\qquad u_t(0,x)=f(x),
\end{equation}
has received a great deal of attention 
(see, among the others,
\cite{B}, \cite{BS}, \cite{BPST}, \cite{C}, \cite{GeV}, \cite{PST}, \cite{Y1}, \cite{Y2} and, for the Schr\"odinger equation, \cite{JSS}, \cite{RS}). Indeed, potential perturbations of the wave equation
are frequently encountered when studying the stability of stationary solutions for several important systems of partial differential equations (wave-Schr\"o\-din\-ger, Max\-well-Schr\"odinger, Maxwell-Dirac and many others). The potentials that arise are usually non-smooth, and this is one of the main motivations for considering rough functions $V(x)$ in \eqref{eq.wavepot}. 

The proof of \eqref{eq.dispfree} is based either on the explicit expression of the fundamental solution, or on the method of stationary phase (see e.g. \cite{SS}). Only partial substitutes of these methods are available
for the perturbed operator. Beals and Strauss \cite{BS}) obtained $L^{p}-L^{p'}$ decay estimates of Strichartz type, later improved by Beals \cite{B} who could prove an almost optimal dispersive estimate similar to \eqref{eq.dispfree}, for smooth positive potentials decaying fast enough at infinity (and $n\geq3$). Their method is based on a repeated use of Duhamel's formula and an explicit representation of the kernels of the operators that arise. The same methods cover
the case of small, smooth, rapidly decaying potentials with undefinite sign. 
For general dimension $n$, the best results are due to Yajima, who, in
a series of papers (see e.g. \cite{Y1}, \cite{Y2}), proved the $L^{p}$ boundedness of the wave operator intertwining the free with the perturbed operator; as a consequence he obtains dispersive estimates for a variety of equations, including the wave equation. We should also mention that
the Strichartz estimates can be proved independently of the dispersive
estimates, under quite general assumptions on the perturbed operator; for a
nice proof see \cite{BPST2}; see also \cite{BPST} and \cite{C}.

In the special case of dimension $n=3$, Georgiev and Visciglia \cite{GeV} were able to prove the dispersive estimate for potentials of H\"older class $V(x)\in C^{\alpha}(\mathbb{R}^{3}\setminus 0)$, $\alpha\in]0,1[$, satisfying for some $\varepsilon>0$
\begin{equation}\label{eq.GV}
    0\leq V(x)\leq\frac C{|x|^{2+\varepsilon}+|x|^{2-\varepsilon}}.
\end{equation}
In particular this implies $V\in L^{3/2-\delta}\cap L^{3/2+\delta}$ for $\delta$ small ($0<\delta<3\varepsilon/4$). 
This decay assumption on the potential $V(x)$ is close to critical, at least in view of the dispersive estimate. Indeed, in \cite{PST}, \cite{PST2} the inverse square potential $V=a|x|^{-2}$ was considered; 
this kind of potential has a critical behaviour, and the dispersive estimate 
\eqref{eq.dispfree} was showed to be false for small negative $a$ (but still true for $a\geq0$). Notice that the inverse square potential belongs to the weak 
$L^{3/2}_{w}\simeq L^{3/2,\infty}$ Lorentz space.

Thus it is natural to ask what are the weakest assumptions on the potential that imply the dispersive estimate. Here we prove that it is sufficient to assume that $V$ belongs to a suitable Kato class of potentials, and no smoothness at all is required. We recall the relevant definitions:

\begin{definition}\label{def.kato}
    The measurable function $V(x)$ on  $\mathbb{R}^{n}$, $n\geq3$, is said to belong to the \emph{Kato class} if
\begin{equation}\label{eq.katoclass}
    \lim_{r\downarrow 0}\sup_{x\in\mathbb{R}^{n}}
      \int_{|x-y|<r}\frac{|V(y)|}{|x-y|^{n-2}}dy=0.
\end{equation}
Moreover, the \emph{Kato norm} of $V(x)$ is defined as
\begin{equation}\label{eq.katonorm}
    \|V\|_{K}=\sup_{x\in\mathbb{R}^{n}}
    \int_{\mathbb{R}^{n}}\frac{|V(y)|}{|x-y|^{n-2}}dy.
\end{equation}
For $n=2$ the kernel $|x-y|^{2-n}$ is replaced by $\log(|x-y|^{-1})$.

The two notions are of course related (e.g., a compactly supported function of Kato class has a finite Kato norm, see Lemma \ref{lem.katoprop} in Section 4).
\end{definition}

\begin{remark}\label{rem.kato}
The relevance of the Kato class in the study of Schr\"odinger operators is well known; full light on its importance was shed in Simon \cite{S} and Aizenmann and Simon \cite{AS}. The stronger norm \eqref{eq.katonorm} was used by Rodnianski and Schlag \cite{RoS} who proved the dispersive estimate for the three dimensional Schr\"odinger equation with a potential having both the Kato and the Rollnik norms small.
\end{remark}

As it is well known, the presence of eigenvalues or resonances
can influence the decay properties of the solutions. The standard
way out of this difficulty is to assume that no resonances are present
on the positive real axis, and in many cases this reduces to assuming
that 0 is not a resonance. In our first result this assumption takes the
following form. We denote as usual by $R_{0}(z)=(-z-\Delta)^{-1}$ 
the resolvent 
operator of $-\Delta$, and by $R_{0}(\lambda\pm i0)$
the limits $\lim_{\varepsilon\downarrow0}R(\lambda\pm i\varepsilon)$
at a point $\lambda\geq0$.
Then we assume that
\begin{quote}
\emph{The integral equation $f+R_{0}(\lambda+i0)Vf=0$
has no nontrivial bounded solution for any $\lambda\geq0$}, 
\end{quote}
or, equivalently,
\begin{equation}\label{eq.reson}
    f+\frac 1{4\pi}
    \int_{\mathbb{R}^3}\frac{e^{i\sqrt\lambda|x-y|}}{|x-y|}V(y)f(y)dy=0,\quad
    f\in L^{\infty},\quad
    \lambda\geq0\quad
    \implies
    \quad
    f\equiv0.
\end{equation}
In several cases this assumption can be
drastically weakened, as discussed below.

We can now state our first result:

\begin{theorem}\label{th.main}
Let $V=V_{1}+V_{2}$ be a real valued potential of Kato class.
Assume that:

i) $V_{1}$ is compactly supported and has a 
bounded Kato norm;

ii) $V_{2}$ has a small Kato norm and precisely
\begin{equation}\label{eq.assv2b}
    \|V_{2}\|_{K}\cdot\left(1+\frac{1}{4\pi}\|V_{1}\|_{K}\right)<4\pi;
\end{equation}

iii) the negative part $V_{-}=\max\{-V,0\}$ satisfies
\begin{equation}\label{eq.negV}
    \|V_{-}\|_{K}<2\pi;
\end{equation}

iv) the non resonant condition \eqref{eq.reson} holds for all 
$\lambda\geq0$.

\noindent
Then any solution $u(t,x)$ to problem \eqref{eq.wavepot} satisfies the dispersive estimate
\begin{equation}\label{eq.free}
   \|u(t,\cdot)\|_{L^{\infty}}\leq 
         C\;t^{-1}\|f\|_{\dot B^{1}_{1,1}(\mathbb{R}^{3})}.
\end{equation}
\end{theorem}

We give some comments on the above assumptions.

\begin{remark}\label{rem.l32}
Condition \eqref{eq.assv2b} can be intepreted as a
smallness at infinity of $V$, and is satisfied by quite a large
class of potentials. For instance,
assume that $V$ belongs to the Lorentz space
$L^{3/2,1}(\mathbb{R}^3)$. By the extended Young inequality 
we have
\begin{equation*}
    \|f\|_{K}\leq c_{0}\|f\|_{L^{3/2,1}}
\end{equation*}
for some universal constant $c_{0}$. Thus we see that
$V$ has a bounded Kato norm,
and a similar argument shows that
$V$ also belongs to the Kato class. Moreover, if
$\chi(x)$ is the characteristic function of the
ball $\{|x|<1\}$,
we can decompose $V$ as follows: for any $R>0$,
\begin{equation*}
    V=V_{1}+V_{2},\quad
    V_{1}=\chi(x/R)V,\quad
    V_{2}=(1-\chi(x/R))V.
\end{equation*}
Notice that
\begin{equation*}
    \|V_{2}\|_{K}\leq c_{0}\|V_{2}\|_{L^{3/2,1}}\to0
     \quad\text{as}\quad R\to+\infty;
\end{equation*}
on the other hand,
\begin{equation*}
    \|V_{1}\|_{K}\leq c_{0}\|V_{1}\|_{L^{3/2,1}}\leq
        c_{0}\|V\|_{L^{3/2,1}}
\end{equation*}
independently of $R$, and hence
\begin{equation*}
        \|V_{2}\|_{K}\cdot\left(1+\frac{1}{4\pi}\|V_{1}\|_{K}\right)\to0
     \quad\text{as}\quad R\to+\infty.
\end{equation*}
In other words, assumptions
(i) and (ii) are automatically satisfied by any potential
in $L^{3/2,1}$. We can sum up this argument in the following
Corollary:

\begin{corollary}\label{cor.main}
    Assume the real valued potential $V$ belongs to $L^{3/2,1}$ with
    $\|V_{-}\|_{K}< 2\pi$ and
    satisfies the non resonant condition \eqref{eq.reson}.
    Then the same conclusion of Theorem \ref{th.main} holds.
\end{corollary}

In particular, this applies to potentials belonging to
$L^{3/2-\delta}(\mathbb{R}^{3})\cap L^{3/2+\delta}(\mathbb{R}^{3})$
for some $\delta>0$, in view of the embedding
\begin{equation*}
     L^{3/2-\delta}(\mathbb{R}^{3})\cap L^{3/2+\delta}(\mathbb{R}^{3})
    \subseteq
    L^{3/2,1}(\mathbb{R}^{3}).
\end{equation*}
This covers the potentials satisfying \eqref{eq.GV}, as remarked above.

It is interesting to compare this to the results of Burq et al. \cite{BPST},
\cite{BPST2} concerning the inverse square potential;
in the scale of Lorentz spaces we can say that the dispersive 
estimate holds when $V\in L^{3/2,1}$ but not when 
$V\in L^{3/2,\infty}$.  It is not clear what can be said for potentials of Lorentz class $L^{3/2,q}$ with $1<q<\infty$, and in particular for
$L^{3/2}=L^{3/2,3/2}$.
\end{remark}

\begin{remark}\label{rem.reson2}
It is a problem of independent interest to find conditions on the potential
$V$ which ensure that no resonances in the
sense of \eqref{eq.reson} occur on the positive real axis.
A well known result in this direction was proved in \cite{ASch} 
(see in particular Appendices 2 and 3). We briefly recall two
special cases which can be applied here ($V$ is always
real valued):

\begin{proposition}\label{prop.ASch} (Alsholm-Schmidt)
Let $n=3$. Assume that $V\in L^{2}_{loc}$ and that, for some $C,R,\epsilon>0$, one has $|V(x)|\leq C{|x|^{-2-\epsilon}}$ for $|x|>R$.
Then property \eqref{eq.reson}
holds for all $\lambda>0$.
\end{proposition}

\begin{proposition}\label{prop.ASch2} (Alsholm-Schmidt)
Let $n=3$. Assume that, for some $C,R,\epsilon>0$, one has $|V(x)|\leq C{|x|^{-1-\epsilon}}$ for $|x|>R$. Moreover, assume that either $V\in L^{1}\cap L^{2}$ 
or $\langle x\rangle^{1/2+\epsilon}V\in L^{2}$.
Then property \eqref{eq.reson}
holds for all $\lambda>0$.
\end{proposition}

Notice that the results of \cite{ASch} do not apply to the potentials
like \eqref{eq.GV} since the singularity $|x|^{-2+\epsilon}$
is not $L^{2}_{loc}$; however, 
in order to apply e.g. Proposition \ref{prop.ASch},
it is sufficient to assume that
\begin{equation}\label{eq.GV2}
    |V(x)|\leq\frac C{|x|^{2+\varepsilon}+|x|^{3/2-\varepsilon}}.
\end{equation}
When $V$ satisfies \eqref{eq.GV2}, (iii) of Theorem \ref{th.main}, and 
$\lambda=0$ is not a resonance (in the sense of \eqref{eq.reson}), then
the dispersive estimate is true.

We further stress that the above propositions do not rule out the possibility of
a resonance at $\lambda=0$. This case can be excluded 
(at least in the sense of \eqref{eq.reson}) if one requires
a stronger decay at infinity of the potential; as an example,
we can prove the following

\begin{theorem}\label{th.mainbis}
Let $V_{1}$ be a nonnegative
$L^{2}$ function such that $V_{1}(x)\leq C|x|^{-3-\delta}$
($\delta>0$)
for large $x$. Then there exists a constant
$\epsilon(V_{1})>0$ such that: for
all real valued functions $V_{2}$ of Kato class with
\begin{equation}\label{eq.katon}
    \|V_{2}\|_{K}<\epsilon(V_{1})
\end{equation}
and for $V=V_{1}+V_{2}$, the solution $u(t,x)$ of problem
\eqref{eq.wavepot} satisfies the dispersive estimate \eqref{eq.free}.
\end{theorem}

In essence, this result states that the dispersive estimate holds
(without additional assumptions on the resonances)
for all nonnegative potentials decaying faster than $|x|^{-3}$
and for all ``small enough'' perturbations thereof; however, it does not
give a measure of the smallness of admissible perturbations. For this, we must
use Theorem \ref{th.main} which requires the additional assumption
\eqref{eq.reson}.

\end{remark}

\begin{remark}\label{rem.paraboliche}
In Section \ref{sec.besov} we prove the equivalence of the standard homogeneous Besov norms with the perturbed ones, i.e., generated by the operator
$-\Delta+V$:
\begin{equation*}
        \dot B^{s}_{1,q}(\mathbb{R}^{n})\cong
        \dot B^{s}_{1,q}(V),\qquad 0<s<2,\quad 1\leq q\leq\infty,\quad
        n\geq3
\end{equation*}
for all potentials $V=V_{+}-V_{-}$ with $V_{\pm}\geq0$ and
\begin{equation}\label{eq.vass}
    \|V_{+}\|_{K}<\infty,\qquad
    \|V_{-}\|_{K}<{\pi^{n/2}}/{\Gamma\left(\frac n 2-1\right)}
\end{equation}
(see Theorem \ref{th.eqhom}). For this result, 
a suitable extension of some lemmas in \cite{JN1}-\cite{JN2} was needed,
which in turn required an improvement in Simon's estimates for the Schr\"odinger semigroup \cite{S}. Indeed, in Proposition \ref{prop.simon} we prove that the semigroup $e^{t(\Delta-V)}$ has an integral kernel $k(t,x,y)$ such that
($n\geq3$)
\begin{equation}\label{eq.ker}
    |k(t,x,y)|\leq
       \frac{(2\pi t)^{-n/2}}{1-2\|V_{-}\|_{K}/c_{n}}e^{-|x-y|^{2}/8t}
\end{equation}
and satisfies the estimate
\begin{equation}\label{eq.ethLpLq}
    \|e^{-tH}\|_{\Lpq}\leq 
        \frac{(2\pi t)^{-\gamma}}{(1-\|V_{-}\|_{K}/c_{n})^{2}},\qquad
        \gamma=\frac n2\left(\frac1p-\frac1q\right).
\end{equation}
Thus, as a byproduct of our proof we obtain the following parabolic dispersive estimate (see Proposition \ref{prop.simon}):

\begin{theorem}\label{th.parab}
Let $n\geq3$, assume the potential $V(x)$ is of Kato class, has a finite Kato norm and its negative part $V_{-}$ satisfies
\begin{equation}\label{eq.vass2}
    \|V_{-}\|_{K}<{2\pi^{n/2}}/{\Gamma\left(\frac n 2-1\right)}
\end{equation}
Then the solution $u(t,x)$ to the perturbed heat equation
\begin{equation}\label{eq.pertheat}
    u_{t}-\Delta u+V(x)u=0,\qquad u(0,x)=f(x)
\end{equation}
satisfies the dispersive estimate
\begin{equation}\label{eq.dispheat}
    \|u(t,\cdot)\|_{L^{q}}\leq C t^{\frac n2\left(\frac 1{q}-
         \frac 1{p}\right)}\|f\|_{L^{p}},\qquad
   \frac 1p+\frac 1q=1,\quad q\in[2,\infty].
\end{equation}
\end{theorem}
       \end{remark}

\begin{remark}\label{rem.vitt}
As noticed in \cite{GeV}, in dimension $n=3$ the spectral representation of the solution and an integration by parts are sufficient to prove the dispersive estimate, provided suitable $L^{1}-L^{\infty}$ estimates for the spectral measure are available. Here we follow a similar line of proof; however,
we prefer to apply the spectral theorem outside the real
axis and to prove estimates which are uniform in the imaginary
part of the parameter. This approach does not
require to extend the limiting absorption principle to the
perturbed operator, as it would be necessary
when working on the real axis. See also the previous work 
\cite{P} where the case of potentials with a {small} 
Kato norm was considered.
\end{remark}

The paper is organized as follows. In Section \ref{sec.free} we recall some basic properties of the free resolvent. Section \ref{sec.pert} is devoted to a study of the operator $-\Delta+V$. In Section \ref{sec.spectral} we prove the crucial estimates for the spectral measure. Section \ref{sec.besov} contains a detailed study of the Schr\"odinger semigroup, which is then applied to estimate functions of the operator $f(-\Delta+V)$, and to prove the equivalence of free and perturbed Besov spaces. Finally, in Section \ref{sec.final} the estimates of Section \ref{sec.spectral} are applied to the spectral representation of the solution, thus obtaining a dispersive estimate in terms of a perturbed Besov norm; in combination with the equivalence result of Section \ref{sec.besov}, this concludes the proof of the Theorems.

\section{Properties of the free resolvent}\label{sec.free}

We start by recalling the well known representation of the free resolvent 
$R_{0}(z)=(-\Delta-z)^{-1}$ in $\mathbb{R}^3$ (see e.g.~\cite{RS}):
\begin{equation}\label{eq.expl}
R_0(\xi^2)g(x)= (-\Delta - \xi^2)^{-1}g=
\begin{cases}
  \displaystyle \frac{1}{4\pi} \int_{\mathbb{R}^3} \frac{e^{i\xi |x-y|}}{|x-y|} g(y) dy
      &\text{ for $\mathop{\mathrm{Im}} \xi > 0$} \\ \\
  \displaystyle \frac{1}{4\pi} \int_{\mathbb{R}^3} \frac{e^{-i\xi |x-y|}}{|x-y|} g(y) dy
      &\text{ for $\mathop{\mathrm{Im}} \xi < 0$.}\\
\end{cases}
\end{equation}
By elementary computations we obtain that for any $\lambda\in\mathbb{R}^{}$ and $\varepsilon>0$
\begin{equation}\label{eq.freepm}
    R_0(\lapm)g(x)=
       \frac{1}{4\pi} \int \frac{e^{\pm i\sqrt{\laep} |x-y|}}{|x-y|} 
           e^{-\varepsilon|x-y|/2\sqrt{\lambda_{\varepsilon}}} g(y) dy
\end{equation}
where
\begin{equation}\label{eq.laep}
    \laep=\frac{\lambda+(\lambda^{2}+\varepsilon^{2})^{1/2}}{2}>0.
\end{equation}
These formulas define bounded operators on $L^{2}$, provided $\varepsilon>0$
or $\lambda<0$. When approaching the positive real axis, i.e., as $\varepsilon\downarrow0$, this property fails; however if we consider the limit operators for $\lambda\geq0$
\begin{equation}\label{eq.idR0}
 R_{0}(\lambda\pm i0)g(x)=
     \frac{1}{4\pi} \int 
          \frac{e^{\pm i\sqrt\lambda |x-y|}}{|x-y|} g(y) dy
\end{equation}
then the \emph{limiting absorption principle} ensures that $R_{0}(\lapmz)$ are bounded from the weighted space 
$L^{2}(\langle x\rangle^{s}dx)$ to $L^{2}(\langle x\rangle^{-s}dx)$ 
for any $s>1$, and actually $R_{0}(\lapm)\to R_{0}(\lapmz)$ 
in the operator norm (see e.g. \cite{A}, \cite{H}). 

For \emph{negative} $\lambda$ the estimates are of course
much stronger since we are in the resolvent set of $-\Delta$. 
Using
\begin{equation*}
    0<\laep<\frac\varepsilon2,\qquad
    \frac\varepsilon{2\sqrt\laep}\geq\sqrt{|\lambda|}\qquad
    \text{for all $\lambda<0$}
\end{equation*}
we have from \eqref{eq.freepm}, for all $\lambda<0$, $\varepsilon\geq0$
\begin{equation}\label{eq.negR0}
    |R_{0}(\lapm)g(x)|\leq 
       \frac{1}{4\pi} \int 
          \frac{e^{-\sqrt{|\lambda|} |x-y|}}{|x-y|} |g(y)| dy
\end{equation}
and actually for $\lambda<0$, $\varepsilon=0$
\begin{equation*}
     R_{0}(\lapmz)g(x)=
       \frac{1}{4\pi} \int 
          \frac{e^{-\sqrt{|\lambda|} |x-y|}}{|x-y|} g(y) dy.   
\end{equation*}

We collect here some immediate consequences of the above representations 
which will be used in the following. Since
\begin{equation}\label{eq.idR0R0ep}
    [R_{0}(\lambda+i\varepsilon)-R_{0}(\lambda-i\varepsilon)]g=
        \frac{i}{2\pi} \int 
          \frac{\sin (\sqrt\laep |x-y|)}{|x-y|}
          e^{-\varepsilon|x-y|/2\sqrt{\lambda_{\varepsilon}}} g(y) dy   
\end{equation}
we can write for all $\lambda\in\mathbb{R}^{}$ and $\varepsilon\geq0$
\begin{equation}\label{eq.estR0R0}
    \|[R_{0}(\lambda+i\varepsilon)-R_{0}(\lambda-i\varepsilon)]g\|_{L^{\infty}}\leq
        \frac{\sqrt\laep}{2\pi} \|g\|_{L^{1}}.
\end{equation}
        
Recalling Definition \ref{def.kato}, a straightforward computation 
shows that
\begin{equation}\label{eq.estR0V}
    \|R_{0}(\lapm)Vg\|_{L^{\infty}}\leq \frac{1}{4\pi}\|V\|_{K}\|g\|_{L^{\infty}}\qquad
                 \forall\lambda\in\mathbb{R}^{},\ \varepsilon\geq0
\end{equation}
for any measurable function $V(x)$, and in a similar way
\begin{equation}\label{eq.estVR0}
    \|VR_{0}(\lapm)g\|_{L^{1}}\leq \frac{1}{4\pi}\|V\|_{K}\|g\|_{L^{1}}\qquad
                 \forall\lambda\in\mathbb{R}^{},\ \varepsilon\geq0.
\end{equation}
Of course for \emph{negative} $\lambda$ we have better estimates:

\begin{lemma}\label{lem.negR0}
    Assume $V$ is of Kato class and has a finite
    Kato norm. Then
    for all $\delta>0$ there exists $C_{\delta}>0$ such that
\begin{equation}\label{eq.estnegR0V}
    \|R_{0}(\lapm)Vg\|_{L^{\infty}}\leq
    \left(\delta+C_{\delta}\frac{\|V\|_{K}}{\sqrt{|\lambda|}}\right)
    \|g\|_{L^{\infty}}\qquad
    \forall\lambda<0,\ \varepsilon\geq0
\end{equation}
and
\begin{equation}\label{eq.estnegVR0}
    \|VR_{0}(\lapm)g\|_{L^{1}}\leq
    \left(\delta+C_{\delta}\frac{\|V\|_{K}}{\sqrt{|\lambda|}}\right)
    \|g\|_{L^{1}}\qquad
    \forall\lambda<0,\ \varepsilon\geq0.
\end{equation}
\end{lemma}

\begin{proof}
By \eqref{eq.negR0} we have 
\begin{equation*}
    |R_{0}(\lapm)Vg(x)|\leq \frac{1}{4\pi}\int
      \frac{|V(y)|}{|x-y|}|g(y)| e^{-\sqrt{|\lambda|}|x-y|}dy.
\end{equation*}
Now for any $r>0$ we can split the integral in two zones $|x-y|<r$ and $\geq r$; for the first piece we have
\begin{equation*}
    \frac{1}{4\pi}\int_{|x-y|<r}
      \frac{|V(y)|}{|x-y|}|g(y)| e^{-\sqrt{|\lambda|}|x-y|}dy\leq
      \frac{1}{4\pi}\int_{|x-y|<r}
      \frac{|V(y)|}{|x-y|} dy\|g\|_{L^{\infty}}
\end{equation*}
and this can be made smaller than $\delta \|g\|_{L^{\infty}}$ by the definition of Kato class \eqref{eq.katoclass}, provided we choose $r<r(\delta)$. With this choice we can estimate the second piece as follows
\begin{equation*}
    \frac{1}{4\pi}\int_{|x-y|\geq r(\delta)}
      \frac{|V(y)|}{|x-y|}|g(y)| e^{-\sqrt{|\lambda|}|x-y|}dy\leq
      \frac{\|g\|_{L^{\infty}}}{4\pi r(\delta)\sqrt{|\lambda|}}\int 
      \frac{|V(y)|}{|x-y|} dy
\end{equation*}
where we have used the inequality $e^{-a}\leq1/ a$, and this proves \eqref{eq.estnegR0V}. Estimate \eqref{eq.estnegVR0} follows by duality.
\end{proof}

We shall also need estimates for the square of the resolvent $R_{0}(\lapm)^{2}$. Since by the resolvent identity
\begin{equation*}
    \frac{\mathrm{d}}{\mathrm{d}z}R_{0}(z)=R_{0}^{2}(z),
\end{equation*}
we have the explicit representations
\begin{equation}\label{eq.idR0ep2}
    R_{0}(\lapm)^{2}g=
      \frac{1}{8\pi}
      \left(\pm\sqrt{\laep}+i\frac{\varepsilon}{2\sqrt{\laep}}\right)^{-1}
      \int e^{\left(\pm i\sqrt{\laep}-\frac{\varepsilon}{2\sqrt{\laep}}\right)|x-y|}g(y)dy
\end{equation}
and\begin{equation}\label{eq.idR02}
    R_{0}(\lapmz)^{2}g=
      \pm\frac{1}{8\pi\sqrt\lambda}
      \int e^{\pm i\sqrt{\lambda}|x-y|}g(y)dy.
\end{equation}
From these relations we obtain immediately the estimate, valid for all $\lambda\in\mathbb{R}^{}$ and $\varepsilon\geq0$ with $(\lambda,\varepsilon)\neq(0,0)$
\begin{equation}\label{eq.estR02}
    \|R_{0}(\lapm)^{2}g\|_{L^{\infty}}\leq
      \frac{1}{8\pi\sqrt\laep}\|g\|_{L^{1}}.
\end{equation}

\section{The perturbed operator}\label{sec.pert}

Properties of Schr\"odinger operators with a potential in the 
Kato class are well known, see e.g. \cite{S}, \cite{HV}, \cite{V}. 
Under the assumptions of the Theorem (and actually 
even under weaker assumptions)
one can prove that $H=-\Delta+V$ defines a self-adjoint 
nonnegative operator on $L^{2}$. For the convenience of the 
reader, we sketch the proof in the following

\begin{lemma}\label{lem.selfadj}
Let $V=V_{+}-V_{-}$ with $V_{\pm}\geq0$ be a measurable function on $\mathbb{R}^3$ satisfying 
\begin{equation}\label{eq.assv1d}
    V_{+}\text{ is of Kato class,}\quad \|V_{-}\|_{K}<4\pi.
\end{equation}
Then the operator $-\Delta+V$ defined on $C^{\infty}_{0}(\mathbb{R}^n)$
extends to a unique nonnegative self-adjoint operator $H=-\Delta +V$ with domain $\mathcal{D}(H)=H^2(\mathbb{R}^3)$ such that
\begin{equation}
\left(\psi, H\psi \right)_{L^2} = \left(\psi,-\Delta\psi \right)_{L^2} + \left(\psi, V\psi \right)_{L^2}\geq0
\qquad\forall \; \psi \in H^2(\mathbb{R}^3).
\end{equation}
\end{lemma}

\begin{proof}
We shall use the KLMN Theorem (see \cite{S}, Vol.II, Theorem 10.17). Thus it is sufficient to verify the following inequality:
\begin{equation}
\label{eq.klmn}
\int_{\mathbb{R}^3} |V(x)| |\varphi(x)|^2 dx \leq a \int_{\mathbb{R}^3} |\nabla \varphi(x)|^2 dx + b \|\varphi\|_{L^2(\mathbb{R}^3)}^2
\end{equation}
for some constants $a<1,\;b\in \mathbb{R}$ and for all test functions $\varphi$ (whence the same inequality is true for all $\varphi\in H^{1}$ which is the domain of the form $-(\Delta\varphi,\varphi)$).

First of all we prove that for some $a\in\mathopen{]}  0,1  \mathclose{[}$ and for \emph{all} $b>0$
\begin{equation}
\label{eq.klmn2}
\int_{\mathbb{R}^3} V_{-}(x) |\varphi(x)|^2 dx
 \leq a \|\nabla \varphi\|^2_{L^2(\mathbb{R}^3)} + b \|\varphi\|_{L^2(\mathbb{R}^3)}^2.
\end{equation}
This is equivalent to
\begin{equation*}
|(V_{-}\varphi,\varphi)_{L^2}| \leq a (\varphi, -\Delta \varphi)_{L^2} + b \|\varphi\|_{L^2}^2 =
a\left\|\left(H_0 + \frac{b}{a}\right)^{\frac{1}{2}} \varphi \right\|_{L^2}^2,
\end{equation*}
where $H_{0}=-\Delta$ is the selfadjoint operator with domain $H^{2}(\mathbb{R}^{3})$.
Thus, writing $g= \left(H_0 + \frac{b}{a}\right)^{\frac{1}{2}} \varphi$, the inequality to be proved takes the form
\begin{equation*}
\left\||V_{-}|^{\frac{1}{2}} \left(H_0 + \frac{b}{a}\right)^{-\frac{1}{2}} g \right\|_{L^2} \leq a \|g\|_{L^2},
\end{equation*}
for some $1>a > 0$ and all $b>0$; and this is equivalent
to prove that
\begin{equation}
\|TT^*\|_{L^2\to L^{2}} = a^2 < 1
\end{equation}
where we introduced the operator
$T = |V_{-}|^{\frac{1}{2}} \left(H_0 + \frac{b}{a}\right)^{-\frac{1}{2}}$ and its adjoint
$$T^* = \left(H_0 + \frac{b}{a}\right)^{-\frac{1}{2}} |V_{-}|^{\frac{1}{2}}.$$

Using the explicit representation
\begin{equation*}
\left(H_0 + \frac{b}{a}\right)^{-1} \varphi = \frac{1}{4 \pi} 
   \int_{\mathbb{R}^3} \frac{ e^{-\sqrt{\frac{b}{a}}|x-y|}}{|x-y|}\varphi(y) dy
\end{equation*}
we can write
\begin{equation*}
\begin{aligned}
\|TT^* \varphi\|_{L^2}^2 =&
\left\||V_{-}|^{\frac{1}{2}} \left(H_0 + \frac{b}{a}\right)^{-1} |V_{-}|^{\frac{1}{2}} \varphi\right\|_{L^2}^2 =\\
=&\frac{1}{(4\pi)^2} \int |V_{-}(x)| \left| \int \frac{ e^{-\sqrt{\frac{b}{a}}|x-y|}}{|x-y|}
|V_{-}(y)|^{\frac{1}{2}} |\varphi(y)| dy \right|^2 dx\\
\end{aligned}
\end{equation*}
and by the Cauchy-Schwartz inequality we have
\begin{equation*}
\leq \frac{1}{(4\pi)^2} \int |V_{-}(x)|
{\left( \int \frac{ e^{-\sqrt{\frac{b}{a}}|x-y|}}{|x-y|} |V_{-}(y)| dy \right)
\left( \int \frac{ e^{-\sqrt{\frac{b}{a}}|x-y|}}{|x-y|} |\varphi(y)|^2 dy \right)} dx.\end{equation*}
Now by definition of Kato norm we have (for all $x$ and any $a,b>0$)
\begin{equation}\label{eq.intkato}
\int \frac{ e^{-\sqrt{\frac{b}{a}}|x-y|}}{|x-y|}|V_{-}(y)|dy\leq
\int\frac{|V_{-}(y)|}{|x-y|}dy\leq\|V_{-}\|_{K}
\end{equation}
which implies
\begin{equation*}
    \|TT^* \varphi\|_{L^2}^2 \leq
    \frac{\|V_{-}\|_{K}}{(4\pi)^2} 
    \int\int |V_{-}(x)|\frac{ e^{-\sqrt{\frac{b}{a}}|x-y|}}{|x-y|} |\varphi(y)|^2 dy dx.
\end{equation*}
Using again \eqref{eq.intkato} we obtain
\begin{equation*}
    \|TT^* \varphi\|_{L^2}^2 \leq
    \frac{\|V_{-}\|_{K}^{2}}{(4\pi)^2} 
    \|\varphi\|_{L^2}^2
\end{equation*}
which means
\begin{equation}\label{eq.tt}
\|TT^*\|_{L^2\to L^{2}} \leq \frac{\|V_{-}\|_K}{4\pi} \equiv a < 1
\end{equation}
by assumption \eqref{eq.assv1d}, and this proves \eqref{eq.klmn2}

To conclude the proof it is sufficient to show that for all test functions $\varphi$, for all $a>0$ and for \emph{some} $b=b(a)\in\mathbb{R}^{}$
\begin{equation}
\label{eq.klmn1}
\int_{\mathbb{R}^3} V_{+}(x) |\varphi(x)|^2 dx \leq a 
\|\nabla\varphi\|_{L^2(\mathbb{R}^3)}^2
+ b \|\varphi\|_{L^2(\mathbb{R}^3)}^2
\end{equation}
The proof is almost identical to the above one; the only difference appears in estimate \eqref{eq.intkato} where we split the integral as follows
\begin{equation*}
    \int \frac{ e^{-\sqrt{\frac{b}{a}}|x-y|}}{|x-y|}|V_{+}(y)|dy
    =\int_{|x-y|<r}+\int_{|x-y|\geq r}
\end{equation*}
for arbitrary $r>0$. Fix now $\delta>0$; if we choose $r>0$ small enough, the first integral can be made smaller than $\delta$ by assumption \eqref{eq.assv1d}; on the other hand, with $r$ chosen, the second integral can be made smaller than $\delta$ by choosing $b$ large enough. In conclusion we have
\begin{equation*}
    \int \frac{ e^{-\sqrt{\frac{b}{a}}|x-y|}}{|x-y|}|V_{+}(y)|dy
    \leq 2 \delta
\end{equation*}
provided $b$ in \eqref{eq.klmn1} is large enough.

Inequality \eqref{eq.klmn} is now a trivial consequence of \eqref{eq.klmn2} and \eqref{eq.klmn1}; thus the assumptions of the KLMN theorem are satisfied and we can construct $H=-\Delta+V$ as a selfadjoint operator on $H^{2}$. To check that it is positive, we write
\begin{equation*}
\left((-\Delta + V)\varphi, \varphi \right)_{L^2} = \left(-\Delta \varphi, \varphi \right)_{L^2} +
\left( V\varphi, \varphi \right)_{L^2} \geq \|\nabla \varphi\|_{L^2}^2 -
 |(V_{-}\varphi,\varphi)_{L^2}|;
\end{equation*}
by inequality (\ref{eq.klmn2}) we may continue
\begin{equation*}
\geq (1-a) \|\nabla \varphi\|_{L^2}^2 - b \|\varphi\|_{L^2}^2 \geq -b \|\varphi\|_{L^2}^2
\end{equation*}
for every $b>0$, and this implies
\begin{equation}
 \left((-\Delta + V)\varphi, \varphi \right)_{L^2}\geq 0.
\end{equation}
\end{proof}

\begin{remark}\label{rem.ngeq3}
The above proof can be easily extended to general dimension $n\geq3$.
Indeed, the kernel $K_{M}(x)$ of $(-\Delta+M)^{-1}$ 
for $M>0$ satisfies
\begin{equation}\label{eq.kerM}
|K(x)|\leq\frac{1}{\alpha_{n}|x|^{n-2}},\qquad
\lim_{M\to+\infty}\sup_{|x|>r}e^{|x|}K(x)=0
\end{equation}
for each fixed $r>0$ (see e.g. \cite{S}, p.454), and
these are exactly the properties we used in the
above proof. Moreover, the
constant $\alpha_{n}$ is well known and is equal to
\begin{equation*}
    \alpha_{n}={4\pi^{n/2}}/\Gamma\left(\frac n2-1\right).
\end{equation*}
Thus we see that the result of Lemma \ref{lem.selfadj}
is true for all $n\geq3$, provided the negative part of $V$
satisfies
\begin{equation}\label{eq.negVself}
    \|V_{-}\|_{K}<{4\pi^{n/2}}/\Gamma\left(\frac n2-1\right).
\end{equation}

\end{remark}

\section{Spectral calculus for the perturbed operator}\label{sec.spectral}

Lemma \ref{lem.selfadj} allows us to apply the spectral theorem and hence to use the functional calculus for $H=-\Delta+V$, i.e., given any function $\phi(\lambda)$ continuous and bounded on $\mathbb{R}^{}$, we can define the operator $\phi(H)$ on $L^{2}$ as
\begin{equation}\label{eq.spectral}
    \phi(H)f=\frac 1{2\pi i}\cdot L^{2}-\lim_{\varepsilon\downarrow0}
       \int\phi(\lambda)[R_{V}(\lap)-R_{V}(\lam)]fd\lambda
\end{equation}
where
\begin{equation*}
    R_{V}(z)=(-\Delta+V-z)^{-1}
\end{equation*}
is the resolvent operator for $H$ (see e.g. Vol. II of \cite{T}). 
When the limit absorption principle is satisfied, one can define the
limit operators $R_{V}(\lapmz)$ and take the limit in the
spectral formula as $\varepsilon\to0$. Instead, here we shall use
formula \eqref{eq.spectral} exclusively, since our estimates
will always be uniform in the parameter $\varepsilon>0$.

For $z$ outside the positive real axis we have the well known identities
\begin{equation}\label{eq.resolvid1}
    R_{0}(z)=\left(I+R_{0}(z)V\right)R_{V}(z)=R_{V}(z)\left(I+VR_{0}(z)\right),
\end{equation}
and a standard way to represent $R_{V}(z)$ in 
terms of $R_{0}(z)$ is to construct the inverse operators
$\left(I+R_{0}(z)V\right)^{-1}$. This is the
content of the following proposition, which is the crucial result of the paper.
In the following we shall consider in detail the case
of dimension 3 alone, but all 
the results in this section can be
extended to general dimension $n\geq2$
by suitable modifications in the proofs.

\begin{proposition}\label{prop.invert}
    Under the assumptions of Theorem \ref{th.main} 
    (or Theorem \ref{th.mainbis}) there exists $\varepsilon_{0}>0$ 
    such that the bounded operators
    $I+R_{0}(\lapm)V\colon L^{\infty}\to L^{\infty}$ are invertible for all 
    $\lambda\in\mathbb{R}^{}$, $0\leq\varepsilon\leq\varepsilon_{0}$ 
    with a uniform bound
\begin{equation}\label{eq.invert}
    \|(I+R_{0}(\lapm)V)^{-1}\|_{\Lii}\leq C
       \text{\ \ for all }\lambda\in\mathbb{R}^{},\
       0\leq\varepsilon\leq\varepsilon_{0}.
\end{equation}
\end{proposition}

We need a few lemmas. First of all we recall the standard
$L^{2}$ weighted estimate of the free resolvent
(see e.g. \cite{A} or Vol.II of \cite{H}; see also \cite{BRV}):

\begin{lemma}\label{lem.agmon}
    For all $\lambda>0$ and $\varepsilon\geq0$, 
    the free resolvent $R_{0}(\lapm)$ is a bounded operator 
    from the weighted $L^{2}(\wx^{2s}dx)$ to the 
    weighted $L^{2}(\wx^{-2s}dx)$ space
    for any $s>1/2$; moreover the following
    estimate holds with a constant
    $C=C(s)$ independent of
    $\varepsilon$, $\lambda$:
\begin{equation}\label{eq.agmon}
    \|\wx^{-s}R_{0}(\lapm)f\|_{L^{2}}\leq\frac{C}{\sqrt\lambda}
        \|\wx^{s}f\|_{L^{2}}.
\end{equation}
\end{lemma}

The following is an elementary but useful
property of Kato class functions:

\begin{lemma}\label{lem.katoprop}
    A compactly supported function of Kato class has a finite Kato norm.
\end{lemma}

\begin{proof}
Let $V(x)$ be of Kato class with support contained in a ball 
$B(0,R)\subseteq \mathbb{R}^3$. Then by definition we have 
the uniform bound
\begin{equation*}
    \int_{|x-y|\leq 1}|V(y)|dy\leq 
    \int_{|x-y|\leq 1}\frac{|V(y)|}{|x-y|}dy\leq C_{0}
\end{equation*}
for some $C_{0}$ independent of $x$; thus, covering the support of $V$
with a finite number of balls of radius 1,
we see that $V\in L^{1}$. Hence we can write
\begin{equation*}
    \int\frac{|V(y)|}{|x-y|}dy\leq
     \int_{|x-y|\leq 1}\frac{|V(y)|}{|x-y|}dy+
     \int_{|x-y|\geq 1}\frac{|V(y)|}{|x-y|}dy\leq C_{0}+\|V\|_{L^{1}}
\end{equation*}
and this concludes the proof.
\end{proof}

The next lemma is sligthly modified from \cite{S}:

\begin{lemma}\label{lem.simon}
    If $V(x)$ is a compactly supported function in the Kato class, 
    then there exists a 
    sequence of functions $\Vep\in C^{\infty}_{0}(\mathbb{R}^{3})$ such that
    $\|\Vep-V\|_{K}\to0$ and $\supp\Vep\downarrow \supp V$ as 
    $\varepsilon\to0$. When $V\geq0$, the functions $V_{\varepsilon}$ can
    be taken nonnegative too.
\end{lemma}

\begin{proof}
By the preceding lemma $V$ has a finite Kato norm, and clearly it
belongs to $L^{1}$. Consider now a sequence of nonnegative radial mollifiers, i.e., let $\rho(x)\in  C^{\infty}_{0}(\mathbb{R}^{3})$ be a nonnegative radial function with support in the ball $\{|x|\leq 1\}$ such that $\int\rho(x)dx=1$, and set $\rep(x)=\varepsilon^{-3}\rho(x/\varepsilon)$. Then we have the following standard properties of the Newton potential $1/|x|$:
\begin{equation}\label{eq.coul1}
    \frac{1}{|x|}*\rep\equiv\frac{1}{|x|}\text{\ \ for }|x|\geq\varepsilon,
\end{equation}
\begin{equation}\label{eq.coul2}
    \frac{1}{|x|}*\rep\leq\frac{1}{|x|}\text{\ \ for all }|x|\neq0.
\end{equation}
Define now $\Vep=V*\rep$; for fixed $x$ we have
\begin{equation*}
    \left|\int\frac{V(y)}{|x-y|}dy-\int\frac{\Vep(z)}{|x-z|}dz\right|=
    \left|\int V(y)\left(\frac{1}{|x-y|}-\int\frac{\rep(z-y)}{|y-z|}dz\right)dy\right|
\end{equation*}
and since by \eqref{eq.coul2} the term in brackets is positive,
\begin{equation*}
    \leq\int |V(y)|\left(\frac{1}{|x-y|}-\int\frac{\rep(z-y)}{|y-z|}dz\right)dy
    \leq \int_{|x-y|< \varepsilon}\frac{|V(y)|}{|x-y|}dy
\end{equation*}
where in the last step we used \eqref{eq.coul1}. Taking the supremum in $x$,
we obtain
\begin{equation*}
    \|V_{\varepsilon}-V\|_{K}\leq \sup_{x\in \mathbb{R}^3}
    \int_{|x-y|< \varepsilon}\frac{|V(y)|}{|x-y|}dy
\end{equation*}
and recalling Definition \ref{eq.katoclass} we conclude that $\|\Vep-V\|_{K}\to0$.
Finally, the support of $\Vep$ is contained in the set of points at 
distance $\leq\varepsilon$ from the support of $V$, and
clearly $V\geq0$ implies $V_{\varepsilon}\geq0$.
\end{proof}

We prove now a property of the squared operator $(R_{0}V)^{2}$:

\begin{lemma}\label{lem.squared}
    Let $V$ be a compactly supported function in the Kato class. Then for all $\lambda>0$,
    $\varepsilon\geq0$ and $\delta>0$ there exists a constant $C_{\delta}$ depending only
    on $\delta$ such that
\begin{equation}\label{eq.squared}
    \|R_{0}(\lapm)VR_{0}(\lapm)Vf\|_{L^{\infty}}\leq
    \left(\delta+\frac{C_{\delta}}{\sqrt{\lambda}}\right)\|f\|_{L^{\infty}}.
\end{equation}
\end{lemma}

\begin{proof}
By the maximum (Phragm\'en-Lindel\"of)
 principle, since $R_{0}(z)$ is holomorphic, it is sufficient to prove the estimate for $\varepsilon=0$, i.e., for the operators $R_{0}(\lapmz)$. If we
approximate $V$ by the sequence of test functions
$\Vep$ constructed in Lemma \ref{lem.simon}, we can write
\begin{equation*}
    R_{0}(\lapmz)VR_{0}(\lapmz)V
    =R_{0}(V-\Vep)R_{0}V+R_{0}\Vep R_{0}(V-\Vep)+R_{0}\Vep R_{0}\Vep
\end{equation*}
and using estimate \eqref{eq.estR0V} we obtain
\begin{equation}\label{eq.first}
    \|R_{0}VR_{0}Vf\|_{L^{\infty}}\leq
    (2\pi)^{-1}\|V\|_{K}\cdot\|V-\Vep\|_{K}\cdot\|f\|_{L^{\infty}}+
         \|R_{0}\Vep R_{0}\Vep f\|_{L^{\infty}}.
\end{equation}
We can choose $\varepsilon=\varepsilon(\delta)$ so small that
\begin{equation*}
    (2\pi)^{-1}\|V\|_{K}\cdot\|V-\Vep\|_{K}
     \leq \frac12\delta,
\end{equation*}
and hence it sufficient to prove \eqref{eq.squared} 
with $V$ replaced by $\Vep$. Now we have
\begin{equation*}
    |R_{0}\Vep R_{0}\Vep f(x)|\leq
     \int_{|x-y|<r}\frac{|\Vep|}{|x-y|}dy\|R_{0}\Vep f\|_{L^{\infty}}+
      \int_{|x-y|\geq r}\frac{|\Vep R_{0}\Vep f|}{|x-y|}dy;
\end{equation*}
the first term clearly satisfies
\begin{equation*}
    \int_{|x-y|<r}\frac{|\Vep|}{|x-y|}dy\leq
     C\int_{|x-y|<r}\frac{dy}{|x-y|}=\sigma(r)\to0
\end{equation*}
since $\Vep$ is bounded, so that we find for all $r>0$
\begin{equation}\label{eq.second}
    |R_{0}\Vep R_{0}\Vep f(x)|\leq
    \sigma(r)\|V\|_{K}\|f\|_{L^{\infty}}+\frac 1 r \|\Vep R_{0} \Vep f\|_{L^{1}}
\end{equation}
where in the last step we used the property
\begin{equation*}
    \int\frac{|\Vep|}{|x-y|}dy\leq\int\frac{|V|}{|x-y|}dy
\end{equation*}
already used in the course of the proof of Lemma \ref{lem.simon}. 
In order to estimate the second term in \eqref{eq.second},
we may write for some $s>1/2$
\begin{equation*}
     \|\Vep R_{0} \Vep f\|_{L^{1}}\leq
     \|\wx^{s}\Vep\|_{L^{2}}\|\wx^{-s}R_{0}\Vep f\|_{L^{2}}
\end{equation*}
and applying Lemma \ref{lem.agmon} we get
\begin{equation*}
     \leq\frac{C}{\sqrt\lambda}
     \|\wx^{s}\Vep\|^{2}_{L^{2}}\|f\|_{L^{\infty}}\leq
     \frac{C_{1}}{\sqrt{\lambda}}\|f\|_{L^{\infty}}
\end{equation*}
since $V_{\varepsilon}$ is in $C^{\infty}_{0}$.
Coming back to \eqref{eq.second}, we obtain
\begin{equation*}
    |R_{0}\Vep R_{0}\Vep f(x)|\leq
    \left(\sigma(r)\|V\|_{K}\|f\|_{L^{\infty}}+
        \frac {C_{1}} r
        \frac{1}{\sqrt\lambda}\right)\|f\|_{L^{\infty}}
\end{equation*}
whence \eqref{eq.squared} follows.
\end{proof}

We prove now a fundamental compactness property:

\begin{lemma}\label{lem.compact}
    Let $V$ be a compactly supported function in the Kato class. Then for all
    $\lambda\in\mathbb{R}^{}$, $\varepsilon\geq0$ the operator 
    $R_{0}(\lapm) V\colon L^{\infty}\to L^{\infty}$  and the operator
    $VR_{0}(\lapm)\colon L^{1}\to L^{1}$ are compact operators. Moreover, if 
    $f\in L^{\infty}$ then the function $R_{0}(\lapm)Vf$ satisfies
    \begin{equation}\label{eq.decay}
    |R_{0}(\lapm)Vf|\leq\frac C{\wx}
    \end{equation}
    for some $C>0$, and hence in particular
    $R_{0}(\lapm)Vf\in L^{2}(\wx^{-2s}dx)$ 
    for all $s>1/2$ and $\lambda,\varepsilon\geq0$.
\end{lemma}

\begin{proof}
If the support of $V$ is contained in the ball $\{|x|\leq M\}$, 
we see that, for all $|x|>2M$ and $y$ in the support of $V$, 
we have $|x-y|\geq|x|-M\geq|x|/2$. Thus by the explicit 
representation of $R_{0}$ we get
\begin{equation*}
    |R_{0}Vf(x)|\leq\int\frac{|V(y)f(y)|}{|x-y|}dy
    \leq\frac 2{|x|}\int|Vf|dy
    \text{\ \ \ for }|x|\geq 2M
\end{equation*}
and recalling that $V\in L^{1}$ we obtain the inequality
\begin{equation}\label{eq.obvineq}
    |R_{0}Vf(x)|\leq\frac 2{|x|}\|V\|_{L^{1}}\|f\|_{L^{\infty}}
    \text{\ \ \ for }|x|\geq 2M,
\end{equation}
From \eqref{eq.obvineq} and the usual estimate
\begin{equation*}
    |R_{0}Vf(x)|\leq\frac{\|V\|_{K}}{4\pi}\|f\|_{L^{\infty}}.
\end{equation*}
we easily deduce the final statement \eqref{eq.decay} and that 
$R_{0}Vf\in L^{2}(\wx^{-2s}dx)$ for all bounded $f$ and $s>1/2$.

In order to prove the compactness property, we may assume
that $V$ is a smooth function with compact support. Indeed, 
by Lemma \ref{lem.simon}, $V$ can be approximated in the 
Kato norm by test functions $\Vep$, so that
$R_{0}V$ is the limit of the sequence of operators
$R_{0}\Vep$ in the $\Lii$ norm, since
\begin{equation*}
    \|R_{0}\Vep-R_{0}V\|_{\Lii}\leq\frac1{4\pi}\|\Vep-V\|_{K}.
\end{equation*}
Thus the compactness of  $R_{0}V$ follows from the compactness
of  $R_{0}\Vep$. A similar argument holds for $VR_{0}$. 
From now on, we shall assume that $V\in C^{\infty}_{0}$.

Let $f_{j}$ be a bounded sequence in $L^{\infty}$; writing
\begin{equation*}
    \nabla_{x}R_{0}Vf(x)=
    \frac 1{4\pi}
    \int V(y)f(y)\nabla_{x}\left(
    \frac{e^{\pm i\sqrt{\laep} |x-y|}}{|x-y|} 
           e^{-\varepsilon|x-y|/2\sqrt{\lambda_{\varepsilon}}}
    \right)dy
\end{equation*}
we immediately obtain a bound for $\|\nabla R_{0}Vf_{j}\|_{L^{\infty}}$,
uniform in $j$ (recall that $V$ now is smooth and compactly supported).
Thus an application of the Ascoli-Arzel\`a 
theorem shows that the sequence $R_{0}Vf_{j}$ is precompact
in the $L^{\infty}$ norm on any bounded set in $\mathbb{R}^{3}$. 
Using this compactness property for small $x$ 
and again inequality \eqref{eq.obvineq} for large $x$, 
by a diagonal procedure we obtain that $R_{0}Vf_{j}$ has a uniformly convergent 
subsequence on the whole $\mathbb{R}^{3}$.

To prove the compactness of $VR_{0}$ we write it as $VR_{0}=A_{r}+B_{r}$ where
\begin{equation}\label{eq.decompA}
    A_{r} g(x)=
       \frac{V(x)}{4\pi} \int \frac{e^{\pm i\sqrt{\laep} |x-y|}}{|x-y|} 
           e^{-\varepsilon|x-y|/2\sqrt{\lambda_{\varepsilon}}}
           \chi_{r}(x-y)g(y) dy
\end{equation}
\begin{equation}\label{eq.decompB}
    B_{r} g(x)=
       \frac{V(x)}{4\pi} \int \frac{e^{\pm i\sqrt{\laep} |x-y|}}{|x-y|} 
           e^{-\varepsilon|x-y|/2\sqrt{\lambda_{\varepsilon}}}
           (1-\chi_{r}(x-y))g(y) dy;
\end{equation}
here $\chi_{r}(y)=\chi(y/r)$ is a cutoff function equal to 1 for $x$ near the origin and vanishing for large $x$. It is easy to show that $B_{r}$ is a compact operator on $L^{1}$; indeed, it is a bounded operator from $L^{1}$ to $W^{1,1}(\Omega)$ for $\Omega$ any bounded open set containing the support of $V$, while $W^{1,1}(\Omega)$ is compactly embedded in $L^{1}(\mathbb{R}^{3})$ by the Rellich-Kondrachov Theorem. Since $\|A_{r}\|_{\Luu}\to0$ as $r\to0$, we regard as above $VR_{0}$ as the uniform limit of compact operators, and this concludes the proof.
\end{proof}

The following version of the same lemma will be useful later on:

\begin{lemma}\label{lem.compact2}
    Assume $V$ satisfies the inequality
    $|V(x)|\leq C\wx^{-3-\delta}$ for some $C,\delta>0$. Then all
    the conclusions of Lemma \ref{lem.compact}
    remain true.
\end{lemma}

\begin{proof}
The estimate follows immediately from the standard inequality
\begin{equation*}
    \int\frac{dy}{\wy^{3+\delta}|x-y|}\leq \frac C{\wx}
\end{equation*}
(see e.g. Appendix 2 of \cite{AS}). The compactness
property is proved as above using the Ascoli-Arzel\`a Theorem.
\end{proof}

We are now ready to prove the main proposition of this section.

\begin{proof}
(of Proposition  \ref{prop.invert}). 
The inversion of $I+R_{0}(z)V:L^{\infty}\to L^{\infty}$ is quite easy
when $\Re z<<0$. Indeed,
Lemma \ref{lem.negR0} 
states that
for all $\delta>0$ there exists a constant $C_{\delta}>0$
such that
\begin{equation*}
    \|R_{0}(\lapm)V\|_{\Lii}\leq \delta
        +C_{\delta}\frac{\|V\|_{K}}{\sqrt{|\lambda|}},\qquad
        \forall\lambda<0,\ \varepsilon\geq0.
\end{equation*}
Hence, in particular,
for $\lambda<-\delta^{2}(C_{\delta}\|V\|_{K})^{-2}$ we have
$\|R_{0}(\lapm)V\|_{\Lii}<2\delta$, and this means that
the norm $\|R_{0}(\lapm)V\|_{\Lii}$ 
tends to 0 for $\lambda\to-\infty$, 
uniformly in $\varepsilon$. Thus $I+R_{0}(\lapm)V$ can be inverted by expansion in Neumann series for any $\varepsilon\geq0$ and any $\lambda<-M$ provided $M>0$ is large enough, and 
the $\Lii$ norm of the inverse operator is bounded by a constant
depending only on $M$ (and $V$).

We now consider the case $\Re z>>0$.
Let $V=V_{1}+V_{2}$ be as in Theorem \ref{th.main},
and write for brevity
\begin{equation*}
    T=R_{0}(z)V_{1},\qquad S=R_{0}(z)V_{2}.
\end{equation*}

We first notice that $I+S$ can be inverted
for all $z\in\mathbb{C}$, with bounded inverse; 
indeed, by \eqref{eq.estR0V} the norm of $S\colon L^{\infty}\to L^{\infty}$ is bounded by $\|V_{2}\|_{K}/(4\pi)$, which is strictly smaller than 1 by assumption \eqref{eq.assv2b}, and the result follows again by a straightforward Neumann series expansion. We thus get for all $z$
\begin{equation}\label{eq.IS}
    \|(I+S)^{-1}\|_{\Lii}\leq \left(1-\|V_{2}\|_{K}/(4\pi)\right)^{-1}.
\end{equation}

We then invert $I+T$ for large $\lambda=\Re z$. Lemma \ref{lem.squared} ensures that $\|T^{2}\|_{\Lii}\to 0$ as $\lambda\to\infty$. This implies that for any $\delta\in]0,1[$ we can find $\lambda_{\delta}$ such that for all 
$\Re z\geq\lambda_{\delta}$, $I-T^{2}$ is invertible with norm
\begin{equation}\label{eq.IT2}
    \|(I-T^{2})^{-1}\|_{\Lii}\leq\frac{1}{1-\delta}.
\end{equation}
Since $I-T$ has norm in $\Lii$ bounded by $1+(4\pi)^{-1}\|V_{1}\|_{K}$ independently of $z$ and
\begin{equation*}
    (I-T)(I-T^{2})^{-1}=(I+T)^{-1},
\end{equation*}
we conclude that also $I+T$ is invertible for any $\Re z\geq\lambda_{\delta}$,
with bound
\begin{equation}\label{eq.IT}
   \|(I+T)^{-1}\|_{\Lii}\leq\frac{1}{1-\delta}(1+\|V_{1}\|_{K}/(4\pi)).
\end{equation}

Consider now for $\Re z\geq\lambda_{\delta}$ the operator
\begin{equation*}
    S(I+T)^{-1};
\end{equation*}
by the usual bound $\|S\|_{\Lii}\leq\|V_{2}\|_{K}/(4\pi)$ and by \eqref{eq.IT} we obtain
\begin{equation*}
    \|S(I+T)^{-1}\|_{\Lii}\leq \frac 1{4\pi}\|V_{2}\|_{K}
        \frac{1}{1-\delta}\left(1+\dfrac{\|V_{1}\|}{4\pi}\right)=\frac\alpha{1-\delta}
\end{equation*}
where the constant $\alpha$, recalling the main assumption \eqref{eq.assv2b}, satisfies
\begin{equation*}
    \alpha\equiv\frac 1{4\pi}\|V_{2}\|_{K}
        \left(1+\dfrac{\|V_{1}\|}{4\pi}\right)<1.
\end{equation*}
Hence we see that
\begin{equation*}
    \|S(I+T)^{-1}\|_{\Lii}\leq \frac\alpha{1-\delta}<1
\end{equation*}
provided $\delta<1-\alpha$, i.e., provided $\lambda_{\delta}$ is large enough. Thus, choosing a value of $\lambda_{\delta}$ large enough, we have that for 
$\Re z\geq\lambda_{\delta}$ the operator
\begin{equation*}
    I+S(I+T)^{-1}
\end{equation*}
is invertible. Finally, writing
\begin{equation*}
    (I+S+T)^{-1}=(I+T)^{-1}(I+S(I+T)^{-1})^{-1},
\end{equation*}
we see that $I+S+T=I+R_{0}V$ is invertible with the bound
\begin{equation}\label{eq.bound1}
    \|(I+R_{0}(z)V)^{-1}\|_{\Lii}\leq\left(1+\dfrac{\|V_{1}\|}{4\pi}\right)
          \frac{1}{1-\alpha-\delta}
\end{equation}
for $\Re z\geq \lambda_{\delta}$.

It remains to invert $I+S+T$ for $-M\leq\Re z\leq\lambda_{\delta}$,
$0\leq\Im z\leq\varepsilon_{0}$
(or $0\geq\Im z\geq-\varepsilon_{0}$),
with a uniform bound. To this end we shall apply Fredholm theory;
notice that the standard analytic Fredholm theory
cannot be applied directly since we are not in the usual Hilbert framework
but we are working in $L^{\infty}$ instead. We proceed in two slightly
different ways according to the set of available assumptions.

\subsection{Case A: assumptions of Theorem \ref{th.main}}

The first step is to prove that $I+S+T:L^{\infty}\to
L^{\infty}$ is injective.
A general argument shows that this is
always the case when $z$ is outside the positive real axis
$[0,+\infty[$, provided $V=V_{1}+V_{2}$ satisfies (i), (ii) of
Theorem \ref{th.main}. To see this, we approximate $V_{1}$ with
a sequence of nonnegative test functions $V_{\delta}$ in such a way that
$\|V_{1}-V_{\delta}\|_{K}\to0$ (see Lemma \ref{lem.simon}); thus we
can decompose $V$ as
\begin{equation*}
    V=V_{\delta}+W_{\delta},\qquad
    0\leq V_{\delta}\in C^{\infty}_{0},\qquad
    \|W_{\delta}\|_{K}=\|V_{2}+V_{1}-V_{\delta}\|_{K}<4\pi
\end{equation*}
for $\delta$ small enough. Assume now that the bounded function
$g$ satisfies the integral equation
\begin{equation*}
    (I+R_{0}(z)V)g=0,\qquad z\not\in \mathbb{R}^+;
\end{equation*}
we shall prove that $g=0$. Indeed, we can rewrite the equation
as follows:
\begin{equation*}
    (I+R_{0}(z)W_{\delta})g=-R_{0}(z)V_{\delta}g\in L^{\infty}.
\end{equation*}
Now, $R_{0}(z)W_{\delta}$ has norm $<1$ as a bounded
operator on $L^{\infty}$, hence we can invert $I+R_{0}(z)W_{\delta}$
and we obtain
\begin{equation*}
    g=-(I+R_{0}(z)W_{\delta})^{-1}R_{0}(z)V_{\delta}g.
\end{equation*}
Note that
\begin{equation*}
    (I+R_{0}(z)W_{\delta})^{-1}R_{0}(z)=(-z-\Delta+W_{\delta})^{-1}
\end{equation*}
is exactly the resolvent operator of $-\Delta+W_{\delta}$, at a point
$z$ outside the spectrum. Moreover, $V_{\delta}g$ is in $L^{2}$, hence
$g=(-z-\Delta+W_{\delta})^{-1}V_{\delta}g$ is in $H^{2}$; since
\begin{equation*}
    (-z-\Delta+V)g=0,\qquad z\not\in\mathbb{R}^+
\end{equation*}
we conclude
that $g\equiv 0$ as claimed.

When $z\in[0,+\infty[$, assumption (iv) of Theorem \ref{th.main} means
exactly that $I+S+T$ is injective on $L^{\infty}$, 
thus we have nothing to prove in this
case, and we obtain that $I+S+T$ is injective for all values of 
$z\in\mathbb{C}$. 

The second step is to prove that $I+S+T$ is invertible.
Recalling that $I+S$ is invertible for all $z$, we can
write
\begin{equation*}
    I+S+T=(I+T(I+S)^{-1})(I+S)
\end{equation*}
which implies that $I+T(I+S)^{-1}$ is also injective for all $z$.
But $T$, and hence $T(I+S)^{-1}$ are compact operators on $L^{\infty}$,
thanks to Lemma \ref{lem.compact}. By Fredholm theory this 
implies that $I+T(I+S)^{-1}$ is invertible, and in conclusion
$I+S+T$ is invertible too and the following identity holds:
\begin{equation}\label{eq.formal}
    (I+S+T)^{-1}=(I+S)^{-1}(I+T(I+S)^{-1})^{-1}.
\end{equation}

The last step is to prove a uniform bound on $(I+S+T)^{-1}$. This
is the content of the following lemma, which is our $L^{\infty}$
replacement for the usual
analytic Fredholm theory in the Hilbert spaces $L^{2}(\wx^{s}dx)$.

\begin{lemma}\label{lem.unifbd}
    Assume $V=V_{1}+V_{2}$, with $V_{1}$
    compactly supported,
    $\|V_{1}\|_{K}<+\infty$, and $\|V_{2}\|_{K}<4\pi$. If the operator
    $I+R_{0}(z)V:L^{\infty}\to L^{\infty}$ is invertible for all $z$ in a
    compact set $D\subset\mathbb{C}^{+}=\{\Re z\geq0\}$
    (or $D\subset\mathbb{C}^{-}$), then
    \begin{equation*}
    \sup_{z\in D}\|(I+R_{0}(z)V)^{-1}\|_{\Lii}<\infty.
    \end{equation*}
\end{lemma}

\begin{proof}
We write as before 
\begin{equation}\label{eq.ST}
    T=R_{0}(z)V_{1},\qquad S=R_{0}(z)V_{2}
\end{equation}
and when
$z_{n}$ is a sequence of points in $\mathbb{C}$ we shall also write
\begin{equation}\label{eq.STn}
    T_{n}=R_{0}(z_{n})V_{1},\qquad S_{n}=R_{0}(z_{n})V_{2}
\end{equation}
Moreover, we shall denote by $L^{\infty}_{K}$ the space of bounded
compactly supported functions, and by $\Lz$ its closure in $L^{\infty}$;
in other words $\Lz$ is the space of bounded functions vanishing at 
infinity, with
the uniform norm.

The proof consists in several steps.

\textsc{Step 1:} $S$ is a bounded operator from $\Lz$ into itself. Indeed, given any
$\phi\in\Lz$, decompose it as
\begin{equation*}
    \phi=\phi_{M}+\psi_{M},\qquad
    \phi_{M}=\phi\cdot \mathbf{1}_{\{|x|<M\}}
\end{equation*}
where $\mathbf{1}_{\{|x|<M\}}$ is the characteristic function of the
ball $\{|x|<M\}$. As in the proof of Lemma \ref{lem.compact}, we have
immediately
\begin{equation}\label{eq.m1}
    |S\phi_{M}(x)|\leq \frac C {|x|}\|V_{2}\|_{L^{1}(|y|\leq M)}
                     \qquad\text{for $|x|>2M$}.
\end{equation}
On the other hand,
\begin{equation}\label{eq.m2}
    \|S\psi_{M}\|_{L^{\infty}}\leq C\|\psi_{M}\|_{L^{\infty}}\to 0
\qquad\text{for $M\to+\infty$}
\end{equation}
since $\phi$ vanishes at infinity. Then, given any $\delta>0$, we may
choose $M=M_{\delta}$ such that $\|\psi_{M}\|_{L^{\infty}}<\delta$;
from \eqref{eq.m1} we obtain
\begin{equation*}
    |S\phi(x)|\leq|S\phi_{M}(x)|+|S\psi_{M}(x)|\leq
         \frac {\|V_{2}\|_{L^{1}}}{|x|}+\delta\qquad
         \text{for $|x|>2M_{\delta}$}
\end{equation*}
and this implies $S\phi\in\Lz$.
         
\textsc{Step 2:} 
If $D\ni z_{n}\to z$ and $\phi\in\Lz$, then $S_{n}\phi\to S\phi$
uniformly on $\mathbb{R}^n$
(with the notations \eqref{eq.STn}). To prove this, 
we notice that
\begin{equation*}
    \frac{|e^{iw_{n}|x-y|}-e^{iw|x-y|}|}{|x-y|}|\leq
       C|w_{n}-w|
\end{equation*}
provided $w_{n},w$ stay in a compact subset of $\mathbb{C}$;
from this, it easily follows that
\begin{equation}\label{eq.exp}
    |(R_{0}(z_{n})-R_{0}(z)f|
     \leq C(D)\cdot |z^{1/2}-z_{n}^{1/2}|\cdot\|f\|_{L^{1}}    
\end{equation}
with the
determination $(\rho e^{i\theta})^{1/2}=\sqrt\rho e^{i\theta/2}$. Now,
let $\phi\in\Lz$; to prove that $S_{n}\phi=R_{0}(z_{n})V_{2}\phi$ 
converges to  $S\phi=R_{0}(z)V_{2}\phi$ uniformly, we decompose
$\phi=\phi_{M}+\psi_{M}$ as in Step 1 and write
\begin{equation*}
    |S_{n}\phi(x)-S\phi(x)|\leq
         |S_{n}\phi_{M}(x)-S\phi_{M}(x)|+|S_{n}\psi_{M}(x)-S\psi_{M}(x)|.
\end{equation*}
The second term is bounded by
\begin{equation*}
    |S_{n}\psi_{M}(x)-S\psi_{M}(x)|\leq \|V_{2}\|_{K}\|\psi_{M}\|_{L^{\infty}}
\end{equation*}
which can be made smaller than $\delta>0$ provided $M>M_{\delta}$,
as in the preceding step. To the first term we apply \eqref{eq.exp} and
we obtain
\begin{equation*}
    |S_{n}\phi_{M}(x)-S\phi_{M}(x)|
         \leq C(D)\cdot |z_{n}^{1/2}-z^{1/2}|\cdot
         \|V_{2}\|_{L^{1}(|y|\leq M)}\|\phi_{M}\|_{L^{\infty}}
\end{equation*}
whence we see that this term tends uniformly to 0 for each fixed $M$,
when $z_{n}\to z$, $z_{n},z\in D$, and this proves the claim.

Note that in Steps 1 and 2 we did not use the assumption 
$\|V_{2}\|_{K}<4\pi$; both properties are true for potentials of arbitrary 
(but bounded) Kato norm; in particular, they hold for $T,T_{n}$.

\textsc{Step 3:} If $D\ni z_{n}\to z$, $\phi\in\Lz$ and $k\geq1$,
then $S^{k}_{n}\phi\to S^{k}\phi$
uniformly on $\mathbb{R}^n$ (where $S^{k}_{n},S^{k}$ are the
$k$-th powers of the operators defined in
\eqref{eq.ST}, \eqref{eq.STn}). It is sufficient
to write
\begin{equation*}
    S_{n}^{k}-S^{k}=\sum_{j=1}^{k}S^{j-1}_{n}(S_{n}-S)S^{k-j}
\end{equation*}
and prove the convergence of each term separately. Indeed,
$S^{k-j}\phi$ is a fixed element of $\Lz$ by Step 1, hence
$(S_{n}-S)S^{k-j}\phi\to0$ uniformly by Step 2, and remarking that
$S_{n}^{j}$ are bounded operators on $L^{\infty}$ with norm
$\|S_{n}^{j}\|\leq \|S_{n}\|^{j}<1$, we conclude that
$S^{j}_{n}(S_{n}-S)S^{k-j}\phi\to0$ uniformly, as claimed.

\textsc{Step 4:} If  $D\ni z_{n}\to z$ and
$\phi\in\Lz$, then $(I+S_{n})^{-1}\phi$ tends
to $(I+S)^{-1}\phi$ uniformly on $\mathbb{R}^n$. To prove this,
note that can write for any $N\geq1$
\begin{equation*}
    (I+S_{n})^{-1}-(I+S)^{-1}=
        \sum_{k=1}^{N}(-1)^{k}(S_{n}^{k}-S^{k})+
        \sum_{k=N+1}^{\infty}(-1)^{k}(S_{n}^{k}-S^{k});
\end{equation*}
the second sum can be estimated in the norm of bounded
operators on $L^{\infty}$ as follows
\begin{equation*}
    \left\|\sum_{k=N+1}^{\infty}(-1)^{k}(S_{n}^{k}-S^{k})\right\|\leq
     \frac{\|S_{n}\|^{N+1}}{1-\|S_{n}\|}+
     \frac{\|S\|^{N+1}}{1-\|S\|}
\end{equation*}
which is smaller than $\delta$ for $N\geq N_{\delta}$ large enough;
on the other hand, we can apply Step 3 to the terms
$S_{n}^{k}-S^{k}$ for $k=1,\dots,N$,
and this concludes the proof of this step.

\textsc{Step 5:} Conclusion of the proof. We know already that
$(I+S)^{-1}$ is well defined with bounded operator norm for all $z$,
hence by the identity
\begin{equation*}
    I+T+S=(I+S)(I+(I+S)^{-1}T)
\end{equation*}
we see that it is sufficient to bound the operator norm of 
$(I+(I+S)^{-1}T)^{-1}$ for $z\in D$. By the uniform
boundedness principle, our claim reduces to the following: given
any sequence $z_{n}$ in $D$, which can be assumed to converge
to $z\in D$, we have that for all $\phi\in L^{\infty}$ there exists
$c(\phi)>0$ such that, for all $n$,
\begin{equation}\label{eq.sup}
    \|(I+(I+S_{n})^{-1}T_{n})^{-1}\phi\|\leq c(\phi)
\end{equation}
(just take any sequence $z_{n}$ such that the norm in \eqref{eq.sup}
converges to the supremum over $D$). We use again the notations
\eqref{eq.ST}, \eqref{eq.STn}.

Indeed, assume by contradiction that there exists $\phi\in L^{\infty}$ such
that
\begin{equation}\label{eq.contr}
    \|(I+(I+S_{n})^{-1}T_{n})^{-1}\phi\|\to\infty\quad\text{as $z_{n}\to z$}
\end{equation}
and consider the renormalized functions
\begin{equation*}
    \psi_{n}=\frac{(I+(I+S_{n})^{-1}T_{n})^{-1}\phi}
            {\|(I+(I+S_{n})^{-1}T_{n})^{-1}\phi\|_{L^{\infty}}}.
\end{equation*}
Clearly we have
\begin{equation}\label{eq.renorm}
    \|\psi_{n}\|_{L^{\infty}}=1,\qquad
    (I+(I+S_{n})^{-1}T_{n})\psi_{n}\to0 \quad\text{in $L^{\infty}$.}
\end{equation}
We have also $\|T_{n}-T\|\to0$, since using again \eqref{eq.exp}
\begin{equation*}
    |(T_{n}-T)\phi|\leq C(D)\cdot
       |z_{n}^{1/2}-z^{1/2}|\cdot\|V_{1}\|_{L^{1}}\|\phi\|_{L^{\infty}}.
\end{equation*}
This and \eqref{eq.renorm} imply
\begin{equation}\label{eq.renorm2}
    \|\psi_{n}\|_{L^{\infty}}=1,\qquad
    (I+(I+S_{n})^{-1}T)\psi_{n}\to0 \quad\text{in $L^{\infty}$.}
\end{equation}
Now, by Lemma \ref{lem.compact}, we know that $T$ is a compact operator
on $L^{\infty}$ and the image of $T$ is contained in $\Lz$ 
(see \eqref{eq.decay}), hence by possibly extracting a subsequence
we obtain that $T\psi_{n}$ converges uniformly to some function
$\zeta\in\Lz$. Now we can write
\begin{equation*}
    (I+S_{n})^{-1}T\psi_{n}=
     (I+S_{n})^{-1}(T\psi_{n}-\zeta) +(I+S_{n})^{-1}\zeta;
\end{equation*}
since $\|(I+S_{n})^{-1}\|<C$ independent
of $n$, the first term converges uniformly to 0, and by Step 4
we obtain that
\begin{equation*}
    (I+S_{n})^{-1}T\psi_{n}\to
     (I+S)^{-1}\zeta
\end{equation*}
uniformly. By \eqref{eq.renorm2}, this implies the uniform
convergence
\begin{equation*}
    \psi_{n}\to -(I+S)^{-1}\zeta=:\psi;
\end{equation*}
notice in particular that $\|\psi\|_{L^{\infty}}=1$.
Summing up, we have proved that
\begin{equation*}
    \psi_{n}\to \psi\equiv-(I+S)^{-1}\zeta,\qquad
    T\psi_{n}\to\zeta\equiv T\psi
\end{equation*}
and this implies 
\begin{equation*}
    \psi+(I+S)^{-1}T\psi=0\quad\text{i.e.}\quad
    (I+S+T)\psi=0
\end{equation*}
which is absurd since $I+T+S$ is invertible and $\|\psi\|_{L^{\infty}}=1$.
\end{proof}

\subsection{Case B: assumptions of Theorem \ref{th.mainbis}}

We note that a potential $V$ satisfying the new assumptions
can be split as $V=V_{1}'+V_{2}'$ with
$V_{1}',V_{2}'$ as in (i), (ii) of Theorem \ref{th.main}
(take $V_{1}'=V$ for $|x|<R$ and 0 outside, with
$R$ large enough).
Thus, for $z\not\in[0,\lambda_{\delta}]$ the same arguments
as in Case A apply; also Lemma \ref{lem.unifbd} can still be
used. Hence it is sufficient to prove that
$I+R_{0}(z)V$ is invertible for $z\in[0,\lambda_{\delta}]$
under the new assumptions.

Since
$V_{1}$ fulfills the conditions of both Propositions
\ref{prop.ASch} and \ref{prop.ASch2}, we see that the operators
$I+R_{0}(\lapmz)V_{1}$ are injective on $L^{\infty}$
for all $\lambda>0$. 

We now prove injectivity also at $\lambda=0$.
Thus, let the bounded function
$f$ satisfy
\begin{equation}\label{eq.eq}
    f(x)+\int\frac{V_{1}(y)f(y)}{|x-y|}dy=0;
\end{equation}
in particular, $f$ is a weak solution of
\begin{equation*}
    \Delta f=V_{1}f\in L^{2}\implies f\in H^{2}.
\end{equation*}
Now, if $V_{1}(x)<C\wx^{-3-\delta}$ for $|x|>M$,
we have immediately, for all $|x|>2M$,
\begin{equation*}
    |f(x)|\leq \|V_{1}\|_{L^{1}(|x|<M)}\|f\|_{L^{\infty}}\frac C{|x|}+
           C\|f\|_{L^{\infty}}\int\frac{dy}{\wy^{3+\delta}|x-y|}
           \leq\frac C{|x|}
\end{equation*}
(see Lemma \ref{lem.compact2} above). 
Differentiating \eqref{eq.eq} we see that 
$\nabla f$ satisfies an analogous integral equation
\begin{equation*}
    \nabla f(x)+\int{V_{1}(y)f(y)} \nabla_{x}\frac1{|x-y|}dy=0
\end{equation*}
which implies
\begin{equation*}
    |\nabla f(x)|\leq C\|f\|_{L^{\infty}}
       \int\frac{|V_{1}(y)|}{|x-y|^{2}}dy.
\end{equation*}
Proceeding as above, we can write for $|x|>2M$
\begin{equation*}
    |\nabla f(x)|\leq \|V_{1}\|_{L^{1}(|x|<M)}\|f\|_{L^{\infty}}\frac C{|x|^{2}}+
           C\|f\|_{L^{\infty}}\int\frac{dy}{\wy^{3+\delta}|x-y|^{2}}
           \leq\frac C{|x|^{2}}
\end{equation*}
thanks to the standard inequality (see \cite{AS})
\begin{equation*}
    \int\frac{dy}{\wy^{3+\delta}|x-y|^{2}}\leq\frac C{\wx^{2}},
\end{equation*}
Thus we have proved that for all $|x|>2M$
\begin{equation}\label{eq.decf}
    |f(x)|\leq\frac{C}{|x|},\qquad
    |\nabla f(x)|\leq \frac{C}{|x|^{2}}.
\end{equation}
Now a standard cutoff trick can be applied
(see the Appendix of \cite{GV}): let $\phi\in C^{\infty}_{0}$
equal to 0 for $|x|>2$ and equal to 1 for $|x|<1$, consider the
identity
\begin{equation*}
    \int\left(|\nabla f|^{2}+V_{1}|f|^{2}\right)
                 \phi\left(\frac y R\right)dy=
       -\frac 1 R\int_{R\leq|y|\leq 2R}
            \nabla\phi\left(\frac y R\right)\cdot\nabla f\cdot\overline fdy
\end{equation*}
and apply the estimates \eqref{eq.decf}
to the right hand member, for $R$ large enough.
We obtain
\begin{equation*}
    \int\left(|\nabla f|^{2}+V_{1}|f|^{2}\right)
                 \phi\left(\frac y R\right)dy
           \leq\frac C{R}
\end{equation*}
and taking the limit as $R\to\infty$ we conclude that $f\equiv0$,
i.e., 0 is not a resonance.

Writing as before $T=R_{0}(z)V_{1}$, we have just proved that
$I+T$ is injective on $L^{\infty}$ for $z\in[0,\lambda_{\delta}]$.
Now we remark that
we can split $V_{1}=V_{1}'+V_{1}''$ as the sum of
a compactly supported function $V_{1}'\in L^{2}$, hence with bounded
Kato norm, and a function $V_{1}''<C\wx^{-3-\delta}$. The corresponding
operators $T=T'+T''$ are compact on $L^{\infty}$ 
by Lemmas \ref{lem.compact},
\ref{lem.compact2} respectively, hence $T$ is compact
and by Fredholm theory we can conclude that $I+T$ is invertible
for all $z\in[0,\lambda_{\delta}]$. Then Lemma
\ref{lem.unifbd} ensures that the operator
norm $(I+T)^{-1}$ is bounded by some constant
$C_{0}$ uniform on $z\in[0,\lambda_{\delta}]$.

Now, writing
\begin{equation*}
    I+T+S=(I+T)(I+(I+T)^{-1}S)
\end{equation*}
we see that in order to invert $I+T+S$ it is sufficient to
invert $I+(I+T)^{-1}S$; since
\begin{equation*}
    \|(I+T)^{-1}S\|\leq \|(I+T)^{-1}\|\cdot \frac{\|V_{2}\|_{K}}{4\pi}
     \leq C_{0}\frac{\|V_{2}\|_{K}}{4\pi}
\end{equation*}
this can be achieved by a Neumann expansion
as soon as the Kato norm of $V_{2}$ is small
enough, i.e.,
\begin{equation*}
    \|V_{2}\|_{K}<\frac{4\pi}{C_{0}}=:\epsilon(V_{1}).
\end{equation*}
This is exactly assumption \eqref{eq.katon}.

Thus we have proved that $I+S+T$ is invertible for all complex $z$,
and a last application of Lemma \ref{lem.unifbd} concludes
the proof of Case B.
\end{proof}

We can now draw some consequences which shall be used
in the following.

\begin{corollary}\label{cor.invert}
     Under the assumptions of Theorem \ref{th.main} 
     (or Theorem \ref{th.mainbis}) there exists 
     $\varepsilon_{0}>0$ such that the
    bounded operators
    $I+VR_{0}(\lapm)\colon L^{1}\to L^{1}$ are invertible for all 
    $\lambda\in\mathbb{R}^{}$, $0\leq\varepsilon\leq\varepsilon_{0}$ with uniform bound
\begin{equation}\label{eq.invert2}
    \|(I+VR_{0}(\lapm))^{-1}\|_{\Luu}\leq 
     C\text{\ \ for all }\lambda\in\mathbb{R}^{},\
     0\leq\varepsilon\leq\varepsilon_{0}.
\end{equation}
\end{corollary}

\begin{proof}
The operators $I+VR_{0}$ are one to one on $L^{1}$ by duality, since by Proposition \ref{prop.invert} the operators $I+R_{0}V$ are onto. They are onto by Fredholm theory, since $VR_{0}$ are compact operators on $L^{1}$ by Lemma \ref{lem.compact}. Finally, the bound on the inverse also follows by duality and the bound \eqref{eq.invert}; indeed, $(L^{1})'=L^{\infty}$ and hence \begin{equation*}
    \|(I+VR_{0})f\|_{L^{1}}=\sup_{\|h\|_{L^{\infty}}=1}\int h (I+VR_{0})f dx=
         \sup_{\|h\|_{L^{\infty}}=1}\int f (I+R_{0}V)h dx.
\end{equation*}
\end{proof}

As a consequence of \eqref{eq.resolvid1} and of 
Proposition \ref{prop.invert}, Corollary \ref{cor.invert} we can write
the standard representation formulas:
\begin{equation}\label{eq.resolvid2}
    R_{V}(z)=(I+R_{0}V)^{-1}R_{0}(z)=R_{0}(z)(I+VR_{0})^{-1}.
\end{equation}
By combining these relations we easily obtain the identity
\begin{multline}\label{eq.RVRV}
    R_{V}(\lap)-R_{V}(\lam)=\\
      = (I+R_{0}(\lam)V)^{-1}(R_{0}(\lap)-R_{0}(\lam))(I+VR_{0}(\lap))^{-1}
\end{multline}
for all $\lambda\in\mathbb{R}^{}$, $\varepsilon\in]0,\varepsilon_{0}]$. Then by the bounds \eqref{eq.estR0R0} and \eqref{eq.invert}, \eqref{eq.invert2} we obtain
\begin{equation}\label{eq.estRVRV}
    \|[R_{V}(\lambda+i\varepsilon)-R_{V}(\lambda-i\varepsilon)]g\|_{L^{\infty}}\leq
       C\sqrt\laep \|g\|_{L^{1}}.
\end{equation}
for all $\lambda\in\mathbb{R}^{}$, $\varepsilon\in]0,\varepsilon_{0}]$. 

Moreover from \eqref{eq.resolvid2} we get
\begin{equation}\label{eq.RV2}
    R_{V}(\lapm)^{2}
    =(I+R_{0}(\lapm)V)^{-1}R_{0}(\lapm)^{2}(I+VR_{0}(\lapm))^{-1}
\end{equation}
and recalling \eqref{eq.estR02} we obtain
\begin{equation}\label{eq.estRV2}
    \|R_{V}(\lapm)^{2}g\|_{L^{\infty}}\leq
      \frac{C}{\sqrt\laep}\|g\|_{L^{1}}    
\end{equation}
for all $\lambda\in\mathbb{R}^{}$, $\varepsilon\in]0,\varepsilon_{0}]$.

\section{Equivalence of Besov norms}\label{sec.besov}

This section is devoted to prove the equivalence of perturbed and 
standard Besov spaces
\begin{equation}\label{eq.eqbesov}
        \dot B^{s}_{1,q}(\mathbb{R}^{3})\cong
        \dot B^{s}_{1,q}(V)
\end{equation}
which holds for $0<s<2$ and $1\leq q\leq \infty$
under our assumptions. An analogous property holds 
also for non homogeneous spaces.

We begin by adapting to our situation a result of Simon \cite{S} (whose proof 
we follow closely). Hoping that estimates \eqref{eq.eth} and \eqref{eq.kernel} 
may be of independent interest, we shall give the proof 
for general dimension $n$. 
If the negative part of the potential is in the Kato class but not small, 
by Theorem B.1.1 of \cite{S} the semigroup is still bounded, but its norm 
may increase exponentially as $t\to\infty$.

\begin{proposition}\label{prop.simon}
    Assume the potential $V=V_{+}-V_{-}$  on $\mathbb{R}^{n}$, $n\geq3$, $V_{\pm}\geq0$, satisfies
\begin{equation}\label{eq.assv3}
    V_{+}\text{\ is of Kato class}
\end{equation}
and
\begin{equation}\label{eq.assv4}
    \|V_{-}\|_{K}<c_{n}\equiv{2\pi^{n/2}}/{\Gamma\left(\frac n 2-1\right)}
\end{equation}
and consider the selfadjoint operator $H=-\Delta+V$. Then for all $t>0$ and $1\leq p\leq q\leq\infty$ the semigroup $e^{-tH}$ is bounded from $L^{p}$ to $L^{q}$ with norm
\begin{equation}\label{eq.eth}
    \|e^{-tH}\|_{\Lpq}\leq 
        \frac{(2\pi t)^{-\gamma}}{(1-\|V_{-}\|_{K}/c_{n})^{2}},\qquad
        \gamma=\frac n2\left(\frac1p-\frac1q\right).
\end{equation}
Moreover, under the stronger assumption
\begin{equation}\label{eq.assv5}
    \|V_{-}\|_{K}<\frac12 c_{n}
\end{equation}
$e^{-tH}$ is an integral operator with kernel $k(t,x,y)$ satisfying
\begin{equation}\label{eq.kernel}
    |k(t,x,y)|\leq
       \frac{(2\pi t)^{-n/2}}{1-2\|V_{-}\|_{K}/c_{n}}e^{-|x-y|^{2}/8t}.
\end{equation}
\end{proposition}

\begin{proof}
In the following we shall use the more convenient notations
\begin{equation}\label{eq.convention}
    H=-\frac12\Delta+V,\qquad H_{0}=-\frac12\Delta;
\end{equation}
thus in the final step it will be necessary to substitute $t\to 2t$ and $V\to V/2$ in order to obtain the correct estimates.

The fundamental tool will be the \emph{Feynman-Ka\v c} formula
\begin{equation}\label{eq.fk}
    (e^{-tH}f)(x)=
        E_{x}\left(\exp\left(-\int_{0}^{t}V(b(s))ds\right)f(b(t))\right)
\end{equation}
which is valid under much more general assumptions (see e.g. \cite{V}). Here $E_{x}$ is the integral over the path space $\Omega$ with respect to the Wiener measure $\mu_{x}$, $x\in \mathbb{R}^{n}$, while $b(t)$ represents a generic path (brownian motion). We shall not need the full power of the theory but only a few basic facts:

i) Given a non negative function $G(x)$ on $\mathbb{R}^{n}$ we have the identity
\begin{equation}\label{eq.reprEx}
    E_{x}\left(\int_{0}^{t}G(b(s))ds\right)=
       \int Q_{t}(x-y)G(y)dy
\end{equation}
where $Q_{t}(x)$ is the function
\begin{equation}\label{eq.Qt}
    Q_{t}(x)=\int_{0}^{t}(2\pi s)^{-n/2}e^{-|x|^{2}/2s}ds.
\end{equation}
It is easy to see by rescaling that
\begin{equation*}
    \int_{0}^{\infty}(2\pi s)^{-n/2}e^{-|x|^{2}/2s}ds=
      \int_{0}^{\infty}\tau^{\frac n 2-2}e^{-\tau}d\tau\frac{|x|^{2-n}}{2\pi^{n/2}}
      =\Gamma\left(\frac n 2-1\right)\frac{|x|^{2-n}}{2\pi^{n/2}}
\end{equation*}
so that by definition of $c_{n}$ (see \eqref{eq.assv4})
\begin{equation}\label{eq.estQt}
    Q_{t}(x)\leq \frac 1{c_{n}|x|^{n-2}}
\end{equation}
and by \eqref{eq.reprEx}
\begin{equation}\label{eq.a}
    E_{x}\left(\int_{0}^{t}G(b(s))ds\right)\leq
       \frac 1{c_{n}}\|G\|_{K}.
\end{equation}
       
ii) {Khasminskii's lemma} (\cite{K}; B.1.2 in \cite{S}): 
\emph{
if $G(x)$ is a non negative function  on $\mathbb{R}^{n}$ such that for some $t$
\begin{equation}\label{eq.kh1}
    \alpha\equiv\sup_{x}E_{x}\left(\int_{0}^{t}G(b(s))ds\right)<1,
\end{equation}
then
\begin{equation}\label{eq.kh2}
    \sup_{x}E_{x}\left(\exp\left(\int_{0}^{t}G(b(s))ds\right)\right)\leq
             \frac 1{1-\alpha}.
\end{equation}
}
An immediate application is the following: if $V_{-}$ satisfies
\begin{equation*}
    \|V_{-}\|_{K}<c_{n}
\end{equation*}
we have
\begin{equation*}
     \alpha\equiv\sup_{x}E_{x}\left(\int_{0}^{t}V_{-}(b(s))ds\right)\leq
              \frac 1{c_{n}}\|V_{-}\|_{K}<1
\end{equation*}
by \eqref{eq.a}, so that
\begin{equation}\label{eq.kh3}
        \sup_{x}E_{x}\left(\exp\left(\int_{0}^{t}V_{-}(b(s))ds\right)\right)\leq
             \frac 1{1-\|V_{-}\|_{K}/c_{n}}.
\end{equation}

These simple facts gives us the first $L^{\infty}-L^{\infty}$ estimate for the semigroup. Indeed, by the Feynman-Ka\v c formula we have
\begin{multline}\label{eq.Lii}
    \|e^{-tH}f\|_{L^{\infty}}=
         \sup_{x\in\mathbb{R}^{n}}
               E_{x}\left(\exp\left(-\int_{0}^{t}V(b(s))ds\right)f(b(t))\right)\leq\\
          \leq\|f\|_{L^{\infty}}
                E_{x}\left(\exp\left(-\int_{0}^{t}|V_{-}(b(s))|ds\right)\right)\leq 
            \frac{\|f\|_{L^{\infty}}}{1-\|V_{-}\|_{K}/c_{n}}.     
\end{multline}

The second step is a $L^{2}-L^{\infty}$ estimate. By the Feynman-Ka\v c formula and the Schwarz inequality
\begin{multline}\label{eq.split}
    |e^{-tH}f(x)|\leq 
       E_{x}\left(\exp\left(-2\int_{0}^{t}V_{-}(b(s))ds\right)\right)^{1/2}
             E_{x}\left(|f(b(t))|\right)^{1/2}\equiv   \\    \equiv
         \left[(e^{-t(H_{0}+2V)}1)(x)\right]^{1/2}\left[e^{-tH_{0}}|f|^{2}\right]^{1/2}
\end{multline}
where in the last step we used again the formula; now $e^{-tH_{0}}$ is the standard heat kernel which has norm $(2\pi t)^{-n/2}$ as an $L^{1}-L^{\infty}$ operator, while we can apply estimate \eqref{eq.Lii} to the operator $e^{-t(H_{0}+2V)}$. We thus obtain
\begin{equation*}
    |e^{-tH}f(x)|\leq 
         \frac{\|1\|_{L^{\infty}}}{1-2\|V_{-}\|_{K}/c_{n}}(2\pi t)^{-n/4}\|f\|_{L^{2}}
\end{equation*}
which implies
\begin{equation}\label{eq.L2i}
    \|e^{-tH}f\|_{L^{\infty}}\leq 
        \frac{(2\pi t)^{-n/4}}{1-2\|V_{-}\|_{K}/c_{n}}\|f\|_{L^{2}},
\end{equation}
provided
\begin{equation*}
    \|V_{-}\|_{K}<\frac{c_{n}}2.
\end{equation*}
By duality, since $e^{-tH}$ is selfadjoint, we obtain the $L^{2}-L^{\infty}$ estimate
\begin{equation}\label{eq.L12}
    \|e^{-tH}f\|_{L^{2}}\leq 
        \frac{(2\pi t)^{-n/4}}{1-2\|V_{-}\|_{K}/c_{n}}\|f\|_{L^{1}};
\end{equation}
using the semigroup property we can write
\begin{equation*}
    e^{-tH}f=e^{-\frac t 2 H}e^{-\frac t 2 H}f
\end{equation*}
and applying \eqref{eq.L2i} first, then \eqref{eq.L12} we obtain
\begin{equation}\label{eq.L1i}
    \|e^{-tH}f\|_{L^{\infty}}\leq 
        \frac{(\pi t)^{-n/2}}{(1-2\|V_{-}\|_{K}/c_{n})^{2}}\|f\|_{L^{1}}.
\end{equation}
Now recalling \eqref{eq.Lii}, by duality and interpolation we obtain
\begin{equation*}
    \|e^{-tH}f\|_{L^{p}}\leq
       \frac{(\pi t)^{-\gamma}}{(1-2\|V_{-}\|_{K}/c_{n})^{2}}\|f\|_{L^{q}}
\end{equation*}
(the constant could be slightly but not essentially improved) with $\gamma$ as in the statement. The change $t\to2t$, $V\to V/2$ gives \eqref{eq.eth}.

Let now $g(x),h(x)$ be bounded functions; the same argument as in \eqref{eq.split} gives
\begin{equation*}
        |e^{-tH}h(x)|\leq 
         \left[(e^{-t(H_{0}+2V)}|h|)(x)\right]^{1/2}\left[e^{-tH_{0}}|h|(x)\right]^{1/2}
\end{equation*}
and multiplying by $g(x)$ and taking the sup we get
\begin{equation}\label{eq.split2}
        \|ge^{-tH}h\|_{L^{\infty}}\leq 
         \|ge^{-t(H_{0}+2V)}|h|\|_{L^{\infty}}^{1/2}
         \|ge^{-tH_{0}}|h|\|_{L^{\infty}}^{1/2}.
\end{equation}
We choose 
\begin{equation*}
    g=\chi_{K_{1}},\ h=f\chi_{K_{2}}
\end{equation*}
where $f(x)$ is a bounded function while $\chi_{K_{1}},\chi_{K_{2}}$ are the characteristic functions of two disjoint compact sets $K_{1},K_{2}$. We may estimate the first factor in \eqref{eq.split2} using \eqref{eq.L1i} as follows
\begin{equation*}
    \|ge^{-t(H_{0}+2V)}|h|\|_{L^{\infty}}\leq
          \|e^{-t(H_{0}+2V)}|h|\|_{L^{\infty}}\leq
          \frac{(\pi t)^{-n/2}}{(1-4\|V_{-}\|_{K}/c_{n})^{2}}\|f\chi_{K_{2}}\|_{L^{1}}
\end{equation*}
while for the second we may use the explicit kernel of $e^{-tH_{0}}$ i.e.,
\begin{equation*}
    (2\pi t)^{-n/2}\exp(-|x-y|^{2}/2t)
\end{equation*}
and we obtain
\begin{equation*}
    \|ge^{-tH_{0}}|h|\|_{L^{\infty}}
         \leq (2\pi t)^{-n/2}\exp(-d^{2}/2t)\|f\chi_{K_{2}}\|_{L^{1}},\qquad
          d=\dist(K_{1},K_{2}).
\end{equation*}
In conclusion we have
\begin{equation}\label{eq.kernel2}
    \|\chi_{K_{1}}e^{-tH}f\chi_{K_{2}}\|_{L^{\infty}}\leq
         \frac{(\pi t)^{-n/2}e^{-d^{2}/4t}}{1-4\|V_{-}\|_{K}/c_{n}}
         \|f\chi_{K_{2}}\|_{L^{1}},
         \qquad
          d=\dist(K_{1},K_{2}).
\end{equation}
By the Dunford-Pettis Theorem (see Tr\`eves \cite{Tr} and A.1.1-A.1.2 in \cite{S}),
this implies at once that $e^{-tH}$ has an integral kernel representation, with kernel
\begin{equation*}
    k(t,x,y)=\frac{(\pi t)^{-n/2}}{1-4\|V_{-}\|_{K}/c_{n}}e^{-|x-y|^{2}/4t}
\end{equation*}
and this concludes the proof (after rescaling back $t\to2t$, $V\to V/2$).
\end{proof}

We shall now use the above kernel representation of the semigroup to improve a result due to Jensen and Nakamura (Theorem 2.1 in \cite{JN1}):

\begin{proposition}\label{prop.jn}
    Assume the Kato class potential $V=V_{+}-V_{-}$  on $\mathbb{R}^{n}$, $n\geq3$, $V_{\pm}\geq0$, satisfies
\begin{equation}\label{eq.assv6}
    \|V_{+}\|_{K}<\infty
\end{equation}
and
\begin{equation}\label{eq.assv7}
    \|V_{-}\|_{K}<\frac12c_{n}\equiv{\pi^{n/2}}/{\Gamma\left(\frac n 2-1\right)}
\end{equation}
and consider the selfadjoint operator $H=-\Delta+V$. Then for any 
$g\in C^{\infty}_{0}(\mathbb{R}^{})$ and any $\theta>0$ the operator $g(\theta H)$ is bounded on $L^{p}(\mathbb{R}^{n})$, $1\leq p\leq\infty$, with norm independent of $\theta$:
\begin{equation}\label{eq.jn}
    \|g(\theta H)\|_{\Lpp} \leq C(p,n,g,V).
\end{equation}
The same property holds for the rescaled operators
\begin{equation}\label{eq.jn2}
    \|g(\Hth)\|_{\Lpp} \leq C(p,n,g,V),
\end{equation}
where $\Hth=-\Delta+\theta V(\sqrt\theta x)$.
\end{proposition}

\begin{proof}
The proof for fixed $\theta$ is contained in \cite{JN2}. In \cite{JN1}, Theorem 2.1, the result was extended to the uniform estimate \eqref{eq.jn} for $0<\theta\leq1$, under assumptions on the potential weaker than ours. Following that proof, in order to extend the result to $\theta\geq1$ it will be sufficient to prove that a few estimates are uniform in $\theta\geq1$. More precisely, consider the rescaled potential
\begin{equation}\label{eq.rescV}
     \Vth(x)=\theta V(\sqrt\theta x);
\end{equation}
notice that the Kato norm is invariant under this transformation:
\begin{equation}\label{eq.resckato}
    \|\Vth\|_{K}\equiv\|V\|_{K}.
\end{equation}
Consider the operator
\begin{equation}\label{eq.hth}
    \Hth=-\Delta+\Vth.
\end{equation}
We proceed exactly as in the proof of
Theorem 2.1 in \cite{JN1}; as remarked there, \eqref{eq.jn} is a
consequence of \eqref{eq.jn2}. Thus we are reduced to
prove that
\begin{equation}\label{eq.hth2}
    \|g(\Hth)\|_{\Lpp} \leq C
\end{equation}
uniformly in $\theta$, and this amounts to prove three estimates
uniformly in $\theta$:

i) a pointwise estimate for the kernel of $e^{-t\Hth}$, 

ii) an $L^{2}-L^{2}$ estimate for the operator $(\Hth+M)^{-1/2}$, $M>0$ a fixed constant (we can take $M=1$ here),

iii) an $L^{2}-L^{2}$ estimate for the operator $\partial_{x}(\Hth+M)^{-1/2}$.

Step i) follows directly from estimate \eqref{eq.kernel}
\begin{equation}\label{eq.kernel3}
    |k_{\theta}(t,x,y)|\leq
       \frac{(2\pi t)^{-n/2}}{1-2\|{\Vth}_{-}\|_{K}/c_{n}}e^{-|x-y|^{2}/4t}.
\end{equation}
which is uniform in $\theta>0$ since by \eqref{eq.rescV}
\begin{equation*}
    \|\Vth_{-}\|_{K}\equiv\|V_{-}\|_{K}
\end{equation*}
does not depend on $\theta$.

Step ii) is trivial since $\|(\Hth+M)^{-1/2}\|_{\Ldd}\leq M^{-1/2}$. To get iii), we must prove that
\begin{equation*}
    \|\partial_{x}(\Hth+M)^{-1/2}f\|_{L^{2}}\leq C\|f\|_{L^{2}}
\end{equation*}
or equivalently
\begin{equation}\label{eq.h2}
    \|g\|_{\dot H^{1}}\leq C\|(\Hth+M)^{1/2}g\|_{L^{2}}
\end{equation}
for some $C$ independent of $\theta>0$. We rewrite \eqref{eq.h2} as
\begin{equation}\label{eq.h3}
      C^{-1}\|g\|_{\dot H^{1}}\leq 
          (-\Delta g,g)+(\Vth g,g)+M\|g\|_{L^{2}}^{2}.
\end{equation}
Clearly \eqref{eq.h3} is implied by
\begin{equation}\label{eq.h4}
      |(\Vth_{-}g,g)|\leq\alpha\|g\|_{\dot H^{1}}+M\|g\|_{L^{2}}^{2},
                  \quad\alpha<1,\quad
       \text{$\alpha$ independent of $\theta$.}
\end{equation}
Now recall \eqref{eq.klmn2}, where we proved the inequality in dimension $n=3$: for all $b>0$
\begin{equation}\label{eq.h5}
    |(V_{2}\varphi,\varphi)|\leq a(-\Delta\varphi,\varphi)+b\|\varphi\|_{L^{2}}
\end{equation}
where by \eqref{eq.tt}
\begin{equation}\label{eq.h6}
    a^{2}=\frac{\|V_{2}\|_{K}}{4\pi}.
\end{equation}
We can now apply \eqref{eq.h5}, \eqref{eq.h6} to $\Vth_{-}$ whose Kato norm is independent of $\theta$:
\begin{equation*}
    a^{2}=\frac{\|\Vth_{-}\|_{K}}{4\pi}=\frac{\|V_{-}\|_{K}}{4\pi}<\frac{c_{3}}{8\pi}
         =\frac14
\end{equation*}
by \eqref{eq.assv7}, and this concludes the proof of iii) in dimension $n=3$. 

The proof for $n\geq3$ is identical; it is sufficient to use again
\eqref{eq.klmn2}, \eqref{eq.tt} which are still true
for general dimension $n$, as noticed in Remark \ref{rem.ngeq3}.
\end{proof}

The following consequence will be useful:

\begin{corollary}\label{cor.jn}
    Assume $V$ satisfies the assumptions of Proposition \ref{prop.jn}, let $\Hth=-\Delta+\theta V(\sqrt\theta x)$, $H_{0}=-\Delta$, and let $\varphi_{j}(s)=\varphi_{0}(2^{-j}s)$, 
$\psi_{j}(s)=\psi_{0}(2^{-j}s)$ be two homogeneous Paley-Littlewood partitions of unity, $j\in\mathbb{Z}$. Then we have the estimates: for all 
$j,k\in\mathbb{Z}$,
\begin{equation}\label{eq.kj}
    \|\varphi_{j}(\sqrt \Hth)\psi_k(\sqrt{H_{0}})\|_{\Luu}\leq
               C 2^{-2j+2k}
\end{equation}
with a constant $C$ independent of $j,k$ and of $\theta>0$. The same estimates hold interchanging $H_{0}$ and $\Hth$.
\end{corollary}

\begin{proof}
We first note two consequences of \eqref{eq.jn}: 
for all $j$, with a constant independent of $j$,
\begin{equation}\label{eq.pp}
    \|\varphi_{j}(\sqrt \Hth)\Hth\|_{\Lpp}\leq C2^{2j},\qquad
    \|\varphi_{j}(\sqrt \Hth)\Hth^{-1}\|_{\Lpp}\leq C2^{-2j}
\end{equation}
and the analogous ones for $H_{0}$ instead of $H$ (indeed, the case $V=0$ is a special case of \eqref{eq.pp}). The first one follows by choosing
\begin{equation*}
    g(s)=\varphi_{0}(\sqrt s)s\quad\Longrightarrow\quad
    g(2^{-2j} \Hth)=\varphi_{j}(\sqrt \Hth)2^{-2j}\Hth;
\end{equation*}
the second one follows by
\begin{equation*}
    g(s)=\varphi_{0}(\sqrt s)s^{-1}\quad\Longrightarrow\quad
    g(2^{-2j} \Hth)=\varphi_{j}(\sqrt \Hth)2^{2j}\Hth^{-1}.
\end{equation*}
Then we can write
\begin{multline*}
    \varphi_{j}(\sqrt \Hth)\psi_k(\sqrt{H_{0}})=
    \varphi_{j}(\sqrt \Hth)\Hth^{-1}\Hth\psi_k(\sqrt{H_{0}})=\\=
    \varphi_{j}(\sqrt \Hth)\Hth^{-1}H_{0}\psi_k(\sqrt{H_{0}})+
        \varphi_{j}(\sqrt \Hth)\Hth^{-1}\Vth\psi_k(\sqrt{H_{0}}).
\end{multline*}
The first term can be estimated immediately using \eqref{eq.pp}:
\begin{equation*}
     \|\varphi_{j}(\sqrt \Hth)\Hth^{-1}H_{0}\psi_k(\sqrt{H_{0}})\|_{\Lpp}\leq
                 C2^{-2j+2k};
\end{equation*}
for the second one we may write
\begin{equation*}
     \|\varphi_{j}(\sqrt \Hth)\Hth^{-1}\Vth\psi_k(\sqrt{H_{0}})\|_{\Lpp}\leq
                 C2^{-2j}\|\Vth\psi_k(\sqrt{H_{0}})\|_{\Lpp}
\end{equation*}
and since
\begin{equation*}
    \Vth\psi_k(\sqrt{H_{0}})=\Vth R_{0}(0)H_{0}\psi_k(\sqrt{H_{0}}),
\end{equation*}
recalling that $\Vth R_{0}$ is a bounded operator on $L^{1}$ (with norm proportional to the Kato norm of $\Vth$ which does not depend on $\theta$) and applying again \eqref{eq.pp} we obtain \eqref{eq.kj}.

For higher dimension $n>3$ the proof is identical; only
in the last step we need the estimate
\begin{equation*}
    \|VR_{0}(0)f\|_{L^{1}}\leq C \|V\|_{K}\|f\|_{L^{1}}
\end{equation*}
which is true for any $n$.
Indeed, $R_{0}(0)$ apart from a constant is the convolution
with the kernel $|x|^{2-n}$, and this gives immediately that $R_{0}(0)V$ is
bounded on $L^{\infty}$ with norm
$C\|V\|_{K}$. By duality we deduce that
$VR_{0}(0)$ is a bounded operator on $L^{1}$ with
the same norm.
\end{proof}

Using Corollary \ref{cor.jn} we can show the equivalence of non homogeneous Besov spaces $B^{s}_{1,q}(V)$ with the standard ones, and later on we shall prove the more delicate result concerning the homogeneous case. We recall the precise definition: given a homogeneous Paley-Littlewood partition of unity $\varphi_{j}(s)=\varphi_{0}(2^{-j}s)$, $j\in\mathbb{Z}$, we set for $p\in[1,\infty]$,
$q\in[1,\infty[$, $s\in\mathbb{R}^{}$
\begin{equation*}
        \|f\|_{\dot B^{s}_{p,q}(V)}=
      \left(
      \sum_{j\in \mathbb{Z}}2^{jsq}\|\varphi_{j}(\sqrt H)f\|^{q}_{L^{p}}
         \right)^{1/q}
\end{equation*}
with obvious modification when $q=\infty$.
On the other hand, if we consider a non homogeneous Paley-Littlewood partition of unity, i.e., $\varphi_{j}$ as above for $j\geq0$, and we set
\begin{equation*}
    \psi_{0}=1-\sum_{j\geq0}\varphi_{j}
\end{equation*}
we have $\psi_{0}\in C^{\infty}_{0}(\mathbb{R}^{n})$,
and we can define the non homogeneous Besov norm as
\begin{equation*}
        \|f\|_{B^{s}_{p,q}(V)}=
      \left(\|\psi_{0}(\sqrt H)f\|^{q}_{L^{p}}+
      \sum_{j\geq0}2^{jsq}\|\varphi_{j}(\sqrt H)f\|^{q}_{L^{p}}
         \right)^{1/q}
\end{equation*}
When $V=0$ we obtain the classical Besov spaces, which we denote simply
by $\dot B^{s}_{p,q}$ and $B^{s}_{p,q}$.

\begin{theorem}\label{th.besov}
    Assume the Kato class potential $V=V_{+}-V_{-}$ on $\mathbb{R}^{n}$, $n\geq3$, $V_{\pm}\geq0$, satisfies
\begin{equation}\label{eq.assv8}
    \|V_{+}\|_{K}<\infty
\end{equation}
and
\begin{equation}\label{eq.assv9}
    \|V_{-}\|_{K}<\frac12c_{n}\equiv{\pi^{n/2}}/{\Gamma\left(\frac n 2-1\right)}
\end{equation}
Then we have the equivalence of norms
\begin{equation}\label{eq.eqbes}
    \|f\|_{ B^{s}_{1,q}(V)}\cong\|f\|_{ B^{s}_{1,q}}
\end{equation}
for all $q\in[1,\infty]$, $0\leq s<2$. 
Moreover, for the rescaled potentials
\begin{equation}\label{eq.rescv}
    V_{\theta}(x)=\theta V(\sqrt\theta x)
\end{equation}
we have the uniform estimates
\begin{equation}\label{eq.unif}
   C^{-1}\|f\|_{ B^{s}_{1,q}}\leq \|f\|_{ B^{s}_{1,q}(V_{\theta})}\leq C\|f\|_{ B^{s}_{1,q}}
\end{equation}
with a constant $C$ independent of $\theta>0$.
\end{theorem}

\begin{remark}\label{rem.bes}
    In order to improve the result and consider higher values of $s\geq2$ stronger smoothness assumptions on the of the potential $V$ are necessary; we shall not pursue this problem here. Also, to prove the equivalence of Besov spaces $B^{s}_{p,q}$ for $p\neq1$, one should prove different bounds for the operator $VR_{0}$ on $L^{p}$; this is possible but quite technical and we limit ourselves to the case $p=1$ which is our main interest here.
\end{remark}

\begin{proof}
We shall limit ourselves to the case $q=1$ and we shall only prove the inequality
\begin{equation}
\|f\|_{{B_{1,1}^s}(\Vth)} \leq C \|f\|_{B_{1,1}^s};
\end{equation}
the proof of the reverse inequality and of the cases $1<q\leq\infty$ are completely analogous. 

In the following we shall drop the index $\theta$ since all the estimates we use (from Proposition \ref{prop.jn} and Corollary \ref{cor.jn}) have constants independent of $\theta>0$.

Using the notations
\begin{equation*}
    D_{V}=\sqrt{H},\quad D=\sqrt{H_{0}}
\end{equation*}
we have
\begin{equation}
  \|f\|_{{B_{1,1}^s}(V)}= \|\psi_0(D_V) f\|_{L^1} +
    \sum_{j=0}^\infty 2^{js} \|\varphi_j(D_V) f\|_{L^1}.
\end{equation}
Using
$$1=\psi_0(D) +\sum _{k\geq0}\varphi_k(D) , $$
we have
\begin{multline*}
\|f\|_{{B_{1,1}^s}(V)}\leq  \|\psi_0(D_V) \psi_0(D) f\|_{L^1} +
 \sum_{k=0}^\infty  \|\psi_0(D_V) \varphi_k(D) f\|_{L^1}  
 +\\+
 \sum_{j=0}^\infty 2^{js} \|\varphi_j(D_V) \psi_0(D) f\|_{L^1} +
 \sum_{j,k \geq 0} 2^{js} \|\varphi_j(D_V) \varphi_k(D) f\|_{L^1} 
 = \\ =
 I+II+III+IV.
\end{multline*}
We estimate separately the four terms.

Since by \eqref{eq.jn2} $\psi_{0}(D_{V})$ is bounded on $L^{1}$, we have for the first term
\begin{equation}
I= \|\psi_0(D_V) \psi_0(D) f\|_{L^1} \leq C \|f\|_{L^1}
\end{equation}
and since
\begin{equation*}
    \|f\|_{L^{1}}\leq\|\psi_0(D) f\|_{L^1}+
                  \sum_{j\geq0}\|\varphi_j(D) f\|_{L^1}
\end{equation*}
this is smaller than $C\|f\|_{B^{s}_{1,1}}$.

The same argument gives for the second term
\begin{equation*}
II= \sum_{k=0}^\infty  \|\psi_0(D_V) \varphi_k(D) f\|_{L^1} \leq
 C \sum_{k=0}^\infty  \|\varphi_k(D) f\|_{L^1} \leq C\|f\|_{B^{s}_{1,1}}
 \end{equation*}

As to the third term, we can write
\begin{equation*}
\sum_{j=0}^\infty 2^{js} \|\varphi_j(D_V) \psi_0(D) f\|_{L^1} =
\sum_{j=0}^\infty 2^{js} \|\varphi_j(D_V)(-\Delta_V)^{-1} (-\Delta_V) \psi_0(D) f\|_{L^1} 
\end{equation*}
and recalling \eqref{eq.pp} used in the proof of the corollary we have (for $s<2$)
\begin{multline*}
    III\leq C\sum_{j \geq 0} 2^{-j(2-s)}\|(-\Delta_{V})\psi_{0}(D)f\|_{L^{1}}=
    C \|(-\Delta_{V})\psi_{0}(D)f\|_{L^{1}}
    \leq \\ \leq
      C\|(-\Delta)\psi_{0}(D)f\|_{L^{1}}+C\|{V}\psi_{0}(D)f\|_{L^{1}}.
\end{multline*}
Now we have
\begin{equation*}
    \|{V}\psi_{0}(D)f\|_{L^{1}}=\|VR_{0}(0)(-\Delta)\psi_{0}(D)f\|_{L^{1}}
          \leq C{\|V\|_{K}}\|(-\Delta)\psi_{0}(D)f\|_{L^{1}}
\end{equation*}
and since $(-\Delta)\psi_{0}(D)$ is bounded in $L^{1}$ by \eqref{eq.jn2}, we conclude that
\begin{equation}
 III \leq
   C_2 \|f\|_{L^1}\leq C_{3}\|f\|_{B^{s}_{1,1}}
\end{equation}
as for the first term.

Finally, we split the fourth term in the two sums for $j\leq k$ and $j>k$:
\begin{equation*}
 IV=\sum_{j,k \geq 0} 2^{js} \|\varphi_j(D_V) \varphi_k(D) f\|_{L^1}=
  \sum_{j\leq k} +\sum_{j> k}. 
\end{equation*}
For $j\leq k$ we use the fact that $\varphi_j(D_V)$ are bounded on $L^{1}$ with uniform norm by \eqref{eq.jn2} and hence
\begin{equation*}
      \sum_{j\leq k} \leq
           C\sum_{k\geq 0}\|\varphi_k(D)  f\|_{L^1}
                     \sum_{0\leq j\leq k}2^{js}=
          2C\sum_{k\geq 0}2^{ks}\|\varphi_k(D)  f\|_{L^1}.
\end{equation*}
For $j>k$, we write
$\varphi_j = \varphi_j( \varphi_{j-1}+ \varphi_j+ \varphi_{j+1})=
 \varphi_j \widetilde{\varphi_j} $
and we have
\begin{equation*}
 \sum_{j>k} 2^{js} \|\varphi_j(D_V) \varphi_k(D) f\|_{L^1}=
  \sum_{j>k} 2^{js} \|\varphi_j(D_V) \varphi_k(D)  
        \widetilde{\varphi_k(D)} f\|_{L^1};
\end{equation*}
now by the corollary we obtain
\begin{equation*}
 \sum_{j>k} 2^{js} \|\varphi_j(D_V) \varphi_k(D) 
             \widetilde{\varphi_k(D)}f\|_{L^1}\leq
 \sum_{j>k}  C 2^{(k-j)(2-s) }2^{ks} \| \widetilde{\varphi_k}f\|_{L^1} 
 \end{equation*}
and since $\sum_{j>k} 2^{(k-j)(2-s)}<1$ we have
\begin{equation}
IV=\sum_{j,k \geq 0} 2^{js} \|\varphi_j(D_V) \varphi_k(D) f\|_{L^1}\leq
C\sum_{k \geq 0} 2^{k} \| \widetilde{\varphi_k(D)}f\|_{L^1}\leq
C\|f\|_{{B_{1,1}^1}(\mathbb{R}^3)}.
\end{equation}
and this concludes the proof.
\end{proof}

We shall finally show that the preceding result implies the equivalence also for homogeneous Besov spaces. Indeed, the uniformity of estimates \eqref{eq.unif} makes it possible to apply a rescaling argument, using the following lemma:

\begin{lemma}\label{lem.resc}
Let $s\in\mathbb{R}^{}$, $p,q,\in[1,\infty]$. The homogeneous $\dot B^{s}_{p,q}(V)$ norm has the following rescaling property with respect to scaling $(S_{\lambda}f)(x)=f(\lambda x)$:
\begin{equation}\label{eq.resc}
    \|S_{\lambda}f\|_{\dot B^{s}_{p,q}(V)}
        =\lambda^{s-\frac n p}\|f\|_{\dot B^{s}_{p,q}(V_{\lambda^{-2}})}
\end{equation}
provided $\lambda=2^{k}$ for some $k\in\mathbb{Z}$.
\end{lemma}

\begin{remark}
A similar property holds also for any positive $\lambda$, with equality replaced by equivalence of norms, however
\eqref{eq.resc} will be sufficient
for our purposes.    
\end{remark}

\begin{proof}
From the identity
\begin{equation*}
    (-\Delta+V(x))S_{\lambda}f(x)=
       \lambda^{2}S_{\lambda}(-\Delta+\lambda^{-2}V(x/\lambda))f(x)
\end{equation*}
we obtain the rule
\begin{equation*}
    \Delta_{V} S_{\lambda}=\lambda^{2}S_{\lambda}\Delta_{V_{\lambda^{-2}}}
\end{equation*}
with the usual notations
\begin{equation*}
     \Delta_{V}=\Delta+V,\qquad V_{\theta}=\theta V(\sqrt\theta x).
\end{equation*}
This implies
\begin{equation*}
        g(-\Delta_{V}) S_{\lambda}
              =S_{\lambda}\ g(-\lambda^{2}\Delta_{V_{\lambda^{-2}}})
\end{equation*}
and in particular for the functions $\phi_{j}(s)=\phi_{0}(2^{-j}s)$, writing as usual $D_{V}=\sqrt{-\Delta_{V}}$,
\begin{equation*}
    \phi_{j}(D_{V})S_{\lambda}=
    \phi_{0}(2^{-j}D_{V})S_{\lambda}=
    S_{\lambda}\phi_{0}(2^{-j}\lambda D_{V_{\lambda^{-2}}}).
\end{equation*}
With the special choice $\lambda=2^{k}$ this can be written
\begin{equation*}
    \phi_{j}(D_{V})S_{2^{k}}=
    S_{2^{k}}\phi_{j-k}(D_{V_{2^{-2k}}}).
\end{equation*}
Hence we have the identity, for $\lambda=2^{k}$,
\begin{equation*}
    \|S_{\lambda}\|_{\dot B^{s}_{p,q}}^{q}=
    \sum_{j\in\mathbb{Z}}2^{jsq}\|\phi_{j}(D_{V})S_{\lambda}f\|^{q}_{L^{p}}=
    \sum_{j\in\mathbb{Z}}2^{jsq}2^{-knq/p}
           \|S_{\lambda}\phi_{j-k}(D_{V_{2^{-2k}}})f\|_{L^{p}}^{q}
\end{equation*}
since $L^{p}$ rescales as $\lambda^{-n/p}$; writing 
$2^{jsq}2^{knq/p}=2^{k(s-n/p)q}2^{(j+k)sq}$ and shifting the sum $j+k\to j$ we conclude the proof.
\end{proof}

Thus we arrive at the final result of this section:

\begin{theorem}\label{th.eqhom}
        Assume the Kato class potential $V=V_{+}-V_{-}$ on $\mathbb{R}^{n}$, $n\geq3$, $V_{\pm}\geq0$, satisfies
\begin{equation}\label{eq.assv10}
    \|V_{+}\|_{K}<\infty
\end{equation}
and
\begin{equation}\label{eq.assv11}
    \|V_{-}\|_{K}<\frac12c_{n}\equiv{\pi^{n/2}}/{\Gamma\left(\frac n 2-1\right)}
\end{equation}
Then we have the equivalence of norms
\begin{equation}\label{eq.eqbeshom}
    \|f\|_{\dot B^{s}_{1,q}(V)}\cong\|f\|_{\dot B^{s}_{1,q}}
\end{equation}
for all $q\in[1,\infty]$, $0< s<2$. 
Moreover, for the rescaled potentials
\begin{equation*}
    V_{\theta}(x)=\theta V(\sqrt\theta x)
\end{equation*}
we have the uniform estimates
\begin{equation}\label{eq.unifhom}
   C^{-1}\|f\|_{\dot B^{s}_{1,q}}\leq \|f\|_{\dot B^{s}_{1,q}(V_{\theta})}
\leq C\|f\|_{\dot B^{s}_{1,q}}
\end{equation}
with a constant $C$ independent of $\theta>0$.
\end{theorem}

\begin{proof}
We shall consider in detail the case $q=1$ only, the remaining cases being completely analogous.

We already know that \eqref{eq.unifhom} holds for dotless Besov spaces. 
Now we need to prove the following inequalities
\begin{equation}\label{eq.nodot}
    C^{-1}\|f\|_{\dot B^{s}_{1,1}(\Vth)}\leq
    \|f\|_{B^{s}_{1,1}(\Vth)}\leq 
       C \|f\|_{\dot B^{s}_{1,1}(\Vth)}+C \|f\|_{\dot B^{0}_{1,1}(\Vth)}
\end{equation}
with a constant $C$ independent of $\theta>0$. 
       
First of all we prove that ($D=\sqrt{-\Delta}$, $D_{{\Vth}}=\sqrt{-\Delta_{{\Vth}}}$)
\begin{equation}\label{eq.1}
    \sum_{j<-1}2^{js}\|\varphi_{j}(D_{{\Vth}})f\|_{L^{1}}\leq
           C \|\psi_{0}(D_{{\Vth}})f\|_{L^{1}}.
\end{equation}
We notice that $\psi_{0}$ is equal to 1 on the support of $\varphi_{j}$ for $j<-1$. Hence $\varphi_{j}=\varphi_{j}\psi_{0}$ for $j<-1$ and we can write
\begin{equation*}
    \|\varphi_{j}(D_{{\Vth}})f\|_{L^{1}}=
    \|\varphi_{j}(D_{{\Vth}})\psi_{0}(D_{{\Vth}})f\|_{L^{1}}\leq C
    \|\psi_{0}(D_{{\Vth}})f\|_{L^{1}}.
\end{equation*}
(we have used the uniform estimates \eqref{eq.jn}-\eqref{eq.jn2}). Thus \eqref{eq.1} follows, provided $s>0$ so that $\sum_{j<-1}2^{js}$ is convergent.

The term for $j=-1$ is estimated in a simple way 
($\varphi_{-1}=\varphi_{-1}(\psi_{0}+\varphi_{1})$)
\begin{multline}\label{eq.2}
    \|\varphi_{-1}(D_{{\Vth}})f\|_{L^{1}}\leq
           \|\varphi_{-1}(D_{{\Vth}})\psi_{0}(D_{{\Vth}})f\|_{L^{1}}
                +\|\varphi_{-1}(D_{{\Vth}})\varphi_{1}(D_{{\Vth}})f\|_{L^{1}}
           \leq \\ \leq
           C \|\psi_{0}(D_{{\Vth}})f\|_{L^{1}}+C\|\varphi_{1}(D_{{\Vth}})f\|_{L^{1}}.
\end{multline}
Clearly, \eqref{eq.1} and \eqref{eq.2} imply immediately the first inequality \eqref{eq.nodot}.

The second inequality in \eqref{eq.nodot} is easier: it is sufficient to prove that
\begin{equation*}
    \|\psi_{0}(D_{{\Vth}})f\|_{L^{1}}
    \leq C \sum_{j\leq 1}\|\varphi_{j}(D_{{\Vth}})f\|_{L^{1}}
\end{equation*}
which follows from $\psi_{0}=\psi_{0}\cdot\sum_{j\leq1}\varphi_{j}$, the triangle inequality, and the boundedness of $\psi_{0}(D_{{\Vth}})$ on $L^{1}$ with uniform norm. This give \eqref{eq.nodot}. Notice that all the constants appearing in the above inequalities are uniform in $\theta>0$.

By \eqref{eq.nodot} and the equivalence \eqref{eq.unif} we can write for $0<s<2$
\begin{equation*}
    \|f\|_{\dot B^{s}_{1,1}}\leq
     C \|f\|_{B^{s}_{1,1}}\leq
     C \|f\|_{B^{s}_{1,1}(\Vth)}\leq
     C\|f\|_{\dot B^{s}_{1,1}(\Vth)}+C\|f\|_{\dot B^{0}_{1,1}(\Vth)}.
\end{equation*}
If we apply this inequality to a rescaled function $S_{2^{k}}f$ and recall Lemma \ref{lem.resc}, we obtain for all $k\in\mathbb{Z}$
\begin{equation*}
    2^{k(s-n)}\|f\|_{\dot B^{s}_{1,1}}\leq
      C2^{k(s-n)}\|f\|_{\dot B^{s}_{1,1}(V_{\theta2^{-2k}})}+
               C2^{-kn}\|f\|_{\dot B^{0}_{1,1}(V_{\theta2^{-2k}})}
\end{equation*}
with constants independent of $k,\theta$;
we can now choose $\theta=2^{2k}\gamma$,
divide by $2^{k(s-n)}$ and let $k\to+\infty$ to obtain
\begin{equation*}
    \|f\|_{\dot B^{s}_{1,1}}\leq
     C\|f\|_{\dot B^{s}_{1,1}(V_{\gamma})}
\end{equation*}
which is the first part of the thesis. The reverse inequality is proved in the same way.
\end{proof}

\section{Conclusion of the proof}\label{sec.final}

By the spectral calculus for $H=-\Delta+V$, given any bounded continuous function $\phi(s)$ on $\mathbb{R}^{}$, we can represent the operator $\phi(H)$ on $L^{2}$ as
\begin{equation}\label{eq.spectral1}
    \phi(H)f=\frac1{2\pi i}\cdot L^{2}-\lim_{\varepsilon\to0}
    \int\phi(\lambda)[R_{V}(\lap)-R_{V}(\lam)]fd\lambda.
\end{equation}
If $\phi=\psi'$ is the derivative of a $C^{1}$ compactly supported function we can integrate by parts obtaining the equivalent form
\begin{equation}\label{eq.spectral2}
    \phi(H)f=\frac i{2\pi}\cdot L^{2}-\lim_{\varepsilon\to0}
     \int\psi(\lambda)[R_{V}(\lap)^{2}-
        R_{V}(\lam)^{2}]fd\lambda.
\end{equation}
        
Now, fix a smooth function $\psi(s)$ with compact support in $]0,+\infty[$ and consider the Cauchy problem
\begin{equation}\label{eq.sist}
\left\{
\begin{array}{l}
    \square u+V(x)u=0,\qquad t\geq0,\ x\in\mathbb{R}^{3}\\
    u(0,t)=0,\qquad 
    u_{t}(0,x)=\psi(H)g\\
\end{array}
\right.
\end{equation}
for some smooth $g$. Then the solution $u$ can be represented as
\begin{equation*}
    u(t,\cdot)=L^{2}-\lim_{\varepsilon\to0}\uep(t,\cdot)
\end{equation*}
where
\begin{equation}\label{eq.repr1}
    \uep(t,x)=\frac1{2\pi i}
    \int_{0}^{\infty}\frac{\sin(t\sqrt\lambda)}{\sqrt\lambda}
    \psi(\lambda)[R_{V}(\lap)-R_{V}(\lam)]gd\lambda
\end{equation}
or equivalently, after integration by parts,
\begin{multline}\label{eq.repr2}
    \uep(t,x)=\frac1{\pi i t}
     \int_{0}^{\infty} \cos(t\sqrt\lambda)
         \psi'(\lambda)[R_{V}(\lap)-R_{V}(\lam)]gd\lambda+\\
         +\frac1{\pi i t}
     \int_{0}^{\infty} \cos(t\sqrt\lambda)
         \psi(\lambda)[R_{V}(\lap)^{2}-R_{V}(\lam)^{2}]gd\lambda. 
\end{multline}
Estimates \eqref{eq.estRVRV} and \eqref{eq.estRV2} applied to \eqref{eq.repr2} give
\begin{equation*}
    \|\uep(t,\cdot)\|_{L^{\infty}}\leq
       \|g\|_{L^{1}}
       \frac C t\int_{0}^{\infty}\left(|\psi'(\lambda)|\sqrt{\laep}
       +\frac{|\psi(\lambda)|}{\sqrt{\laep}}\right)d\lambda
\end{equation*}
and recalling that
\begin{equation*}
    \lambda\leq\laep\leq\lambda+\frac{\varepsilon}2
\end{equation*}
we obtain
\begin{equation}\label{eq.basic}
    \|\uep(t,\cdot)\|_{L^{\infty}}\leq
       \|g\|_{L^{1}}
       \frac C t\int_{0}^{\infty}\left(|\psi'(\lambda)|(\sqrt{\lambda}+\sqrt\varepsilon)
       +\frac{|\psi(\sqrt\lambda)|}{\sqrt{\lambda}}\right)d\lambda.
\end{equation}

Let now $\varphi_{j}(s)$, $j\in\mathbb{Z}$ be the homogeneous Paley-Littlewood partition of unity defined in the Introduction, with
\begin{equation*}
    \varphi_{j}(s)=\phi_{0}(2^{-j}s),
\end{equation*}
define
\begin{equation}\label{eq.tilde}
    \widetilde\varphi_{j}(s)=\varphi_{j-1}(s)
            +\varphi_{j}(s)+\varphi_{j+1}(s)
\end{equation}
and choose in \eqref{eq.sist}
\begin{equation*}
    \psi(\lambda)=\widetilde\varphi_{j}(\sqrt\lambda)
      \equiv\widetilde\varphi_{0}(2^{-j}\sqrt\lambda).
\end{equation*}
We thus obtain
\begin{equation*}
    \|\uep(t,\cdot)\|_{L^{\infty}}\leq
       \|g\|_{L^{1}}
       \frac C t\int_{0}^{\infty}
       \left(2^{-j}|\widetilde\varphi_{0}'(2^{-j}\sqrt\lambda)|
       \frac{\sqrt{\lambda}+\sqrt\varepsilon}{2\sqrt\lambda}
       +\frac{|\widetilde\varphi_{0}(2^{-j}\sqrt\lambda)|}{\sqrt{\lambda}}\right)
       d\lambda
\end{equation*}
which after the change of variables $\mu=2^{-j}\sqrt\lambda$ gives
\begin{equation}\label{eq.basic2}
    \|\uep(t,\cdot)\|_{L^{\infty}}\leq
          \frac C t (2^{j}+\sqrt\varepsilon)\|g\|_{L^{1}}.
\end{equation}
for some constant $C$ independent of $j,t$ and $g$.
If we let $\varepsilon\to0$, for fixed $t$ the functions $\uep(t,\cdot)$ converge in $L^{2}$ to the solution $u(t,x)$; hence a subsequence converges a.e. and we obtain the estimate
\begin{equation}\label{eq.basic3}
    \|u(t,\cdot)\|_{L^{\infty}}\leq C
          \frac {2^{j}} t \|g\|_{L^{1}}
\end{equation}
for the solution $u(t,x)$ of the Cauchy problem
\begin{equation}\label{eq.sistj}
\left\{
\begin{array}{l}
    \square u+V(x)u=0,\qquad t\geq0,\ x\in\mathbb{R}^{3}\\
    u(0,t)=0,\qquad 
    u_{t}(0,x)=\widetilde\varphi_{j}(\sqrt H)g\\
\end{array}
\right.
\end{equation}
If we now choose
\begin{equation*}
    g=\varphi_{j}(\sqrt H)f
\end{equation*}
and notice that 
$\widetilde\varphi_{j}g\equiv\widetilde\varphi_{j}\varphi_{j}f\equiv\varphi_{j}f$
since $\widetilde\varphi_{j}=1$ on the support of $\varphi_{j}$, we conclude that: the solution $u(t,x)$ of the Cauchy problem
\begin{equation}\label{eq.sisttrue}
\left\{
\begin{array}{l}
    \square u+V(x)u=0,\qquad t\geq0,\ x\in\mathbb{R}^{3}\\
    u(0,t)=0,\qquad 
    u_{t}(0,x)=\varphi_{j}(\sqrt H)f\\
\end{array}
\right.
\end{equation}
satisfies the estimate
\begin{equation}\label{eq.basic4}
    \|u(t,\cdot)\|_{L^{\infty}}\leq C
          \frac {2^{j}} t \|\varphi_{j}(\sqrt H)f\|_{L^{1}}
\end{equation}
          
Consider now the original Cauchy problem \eqref{eq.wavepot}; decomposing the initial datum $f$ as
\begin{equation*}
    f=\sum_{j\in\mathbb{Z}}\varphi_{j}(\sqrt H)f
\end{equation*}
applying estimate \eqref{eq.basic4} and summing over $j$, we obtain by linearity that the solution $u(t,x)$ to \eqref{eq.wavepot} satisfies the estimate
\begin{equation}\label{eq.basic5}
     \|u(t,\cdot)\|_{L^{\infty}}\leq 
          \frac {C} t \|f\|_{\dot B^{1}_{1,1}(V)}.
\end{equation}
Since by Theorem \ref{th.eqhom} this norm is equivalent to the standard one, we see that the proof of Theorem \ref{th.main} is concluded.




\end{document}